\documentclass[11pt,a4paper]{article}
\usepackage[T1]{fontenc}
\usepackage{lmodern,microtype}
\usepackage[margin=27mm]{geometry}
\usepackage{amsmath,amssymb,amsthm,mathtools}
\usepackage{enumitem,needspace}
\usepackage[hidelinks]{hyperref}

\newtheorem{theorem}{Theorem}
\newtheorem{lemma}{Lemma}[section]

\theoremstyle{definition}

\numberwithin{equation}{section}
\newcommand{\Z}{\mathbb Z}
\newcommand{\Q}{\mathbb Q}
\newcommand{\C}{\mathbb C}
\newcommand{\Zp}{\Z_{(p)}}
\newcommand{\Rp}{\mathcal R}
\newcommand{\supp}{\operatorname{supp}}

\newcommand{\ord}{\operatorname{ord}}
\newcommand{\diag}{\operatorname{diag}}
\newcommand{\one}{\mathbf 1}
\setlist[enumerate]{leftmargin=*,topsep=4pt,itemsep=4pt}
\allowdisplaybreaks[1]
\renewcommand{\footnoterule}{\kern-3pt\hrule width 5pc\kern2.6pt}
\hypersetup{pdftitle={Fuglede's conjecture holds on Z\_(p\textasciicircum a q\textasciicircum b)},
  pdfauthor={Zi-Jian Song and Wei-Jie Wang},
  pdfsubject={Spectrality and translational tiling in finite cyclic groups}}

\begin{document}
\begin{center}
{\LARGE Fuglede's conjecture holds on $\Z_{p^a q^b}$\par}
\bigskip
{\footnotesize\MakeUppercase{Zi-Jian Song and Wei-Jie Wang}\par}
\end{center}
\vspace{0.5em}

\begin{abstract}
Let $p,q$ be distinct primes and let $a,b$ be positive integers. We prove
that every nonempty spectral subset of $\Z_{p^a q^b}$ tiles the group. Together with
the known tiling-to-spectral implication, this establishes Fuglede's
conjecture on $\Z_{p^a q^b}$. We work within Malikiosis's reduction to a
maximal spectral pair at a minimal counterexample modulus. Weighted
character projections give a denominator restriction on residue
proportions, and the cyclotomic zeros give a second restriction on
their differences. These restrictions force uniformity modulo $q$
through the initial $p$-levels, in contradiction with absorption-freeness.
\end{abstract}

\begingroup
\makeatletter
\renewcommand{\@makefntext}[1]{\parindent=1em\indent#1}
\makeatother
\renewcommand{\thefootnote}{}
\footnotetext[0]{2020 \textit{Mathematics Subject Classification.}\enspace
  Primary 43A46. Secondary 05B45.\par
  \textit{Key words and phrases.}\enspace
  Fuglede's conjecture, spectral sets,
  Hadamard matrix, weighted character projections,
  cyclotomic zeros.}
\endgroup

\tableofcontents
\clearpage

\section*{Introduction}
\addcontentsline{toc}{section}{Introduction}
Let $\Omega\subseteq\mathbb R^d$ be a bounded measurable set of positive
Lebesgue measure $|\Omega|$. We call $\Omega$ \emph{spectral} if
$L^2(\Omega)$ has an orthonormal basis of the form
\[
 \bigl\{|\Omega|^{-1/2}e^{2\pi i\lambda\cdot x}:\lambda\in\Lambda\bigr\},
\]
where the countable set $\Lambda\subseteq\mathbb R^d$ is called a
\emph{spectrum} of $\Omega$. We say that $\Omega$ \emph{tiles}
$\mathbb R^d$ by translations if there is a countable set
$\mathcal T\subseteq\mathbb R^d$ such that
$\sum_{t\in\mathcal T}\one_\Omega(x-t)=1$ for almost every
$x\in\mathbb R^d$. Fuglede~\cite{Fuglede}
proved that $\Omega$ tiles by a full-rank lattice $L$ if and only if
the dual lattice of $L$ is a spectrum of $\Omega$. He conjectured that
$\Omega$ is spectral if and only if it tiles $\mathbb R^d$ by translations.

For general bounded measurable sets, both directions of the conjecture
fail in every dimension $d\ge2$ (see
\cite{FarkasMatolcsiMora,FarkasRevesz,KolountzakisMatolcsiHadamard,KolountzakisMatolcsiTiles,Matolcsi4D,Tao}
for $d\ge3$ and the preprint~\cite{ZhangPlane} for $d=2$).
For convex bodies the conjecture holds in every
dimension~\cite{GreenfeldLev,IosevichKatzTao,LevMatolcsi}.
On the real line it remains open in general, although it holds for
unions of two or three bounded intervals~\cite{BoseThreeIntervals,LabaTwoIntervals,ShiThreeIntervals}.

The connection with finite cyclic groups comes from the structure of
spectra on the line. Let $\Omega\subset\mathbb R$ be a bounded spectral
set of measure one, and translate a spectrum $\Lambda$ of $\Omega$ so
that $0\in\Lambda$. Then $\Lambda$ has a positive integer period ~\cite{BoseMadan,IosevichKolountzakis}. Dutkay and Lai~\cite[Theorem~1.3]{DutkayLai} proved that, in the
tiling-to-spectral direction, validity on $\mathbb R$ is equivalent to
validity for finite sets in $\Z$ and to validity in all finite cyclic
groups $\Z_N$. They obtained the same equivalence in the
spectral-to-tiling direction under the hypothesis that every bounded
spectral set of measure one on $\mathbb R$ admits a rational spectrum.
The rationality theorem of Fu and Song~\cite{FuSong} supplies this
hypothesis, completing the equivalence in both directions of Fuglede's
conjecture.

\Needspace{7\baselineskip}
Write $\Z_N=\Z/N\Z$ and $\zeta_N=e^{2\pi i/N}$.
For nonempty $A,B\subseteq\Z_N$, we call $(A,B)$ a \emph{spectral pair}
if $|A|=|B|$ and the normalized character matrix
\[
U=|A|^{-1/2}(\zeta_N^{xy})_{x\in A,\,y\in B}
\]
is unitary. In this case $A$ is spectral with spectrum
$B$. We say that $A$ \emph{tiles} $\Z_N$ if there is a set
$T\subseteq\Z_N$ such that each element of $\Z_N$ has a unique
representation $a+t$ with $a\in A$ and $t\in T$. We then write
$\Z_N=A\oplus T$.

For $N=p^n$, the tiling-to-spectral implication follows from Lagarias
and Wang's work on universal spectra~\cite{LagariasWangFactorizations}.
The converse follows from {\L}aba's theorem for spectra with a common
prime-power denominator~\cite[Theorem~1.5(ii)]{Laba} and the
Coven--Meyerowitz tiling criterion~\cite{CovenMeyerowitz}.

Shi~\cite{ShiPQR} proved the conjecture for $N=pqr$, and Kiss et
al.~\cite{KissPQRS} proved it for $N=pqrs$. Somlai~\cite{SomlaiP2QR}
proved the spectral-to-tiling implication for $N=p^2qr$, and
Zhang~\cite{ZhangGroupRing} established it for $N=p^nqr$.
In preprints, Fallon, Kiss, Mayeli, and Somlai~\cite{FallonP2Q2R}
prove the spectral-to-tiling implication for $N=p^2q^2r$ when
$r>p^2q^2$, and Tan and Zhang~\cite{TanZhangSquareFree} report a proof
for all square-free orders. The primes denoted by different letters
are distinct, and $n\geq1$.

\subsection*{From tiling to spectrality}

Coven and Meyerowitz~\cite{CovenMeyerowitz} gave two cyclotomic
conditions, \textup{(T1)} and \textup{(T2)}, which characterize finite
integer tiles whose cardinalities have at most two distinct prime
divisors. {\L}aba~\cite{Laba} showed that sets satisfying these
conditions are spectral. In finite cyclic groups this gives the
tiling-to-spectral implication for orders with at most two distinct
prime divisors. Further cases
include orders $p_1^n p_2\cdots p_r$~\cite[Theorem~1.5]{Malikiosis}
and $(p_1p_2p_3)^2$~\cite{LabaLondnerOdd,LabaLondnerEven},
with the $p_i$ distinct.

In two recent preprints, Li et al.~\cite{LiIntegerTiles} and
Prajapati~\cite{PrajapatiT2} present proofs of the Coven--Meyerowitz
conjecture, which asserts that \textup{(T2)} is necessary for every
finite integer tile. The resulting characterization yields the
tiling-to-spectral implication for every finite cyclic group.

\subsection*{The spectral-to-tiling direction}

Malikiosis and Kolountzakis~\cite{MalikiosisKolountzakis} proved the
spectral-to-tiling implication for $N=p^a q$ using the structure of
vanishing sums of roots of unity. Kiss, Malikiosis, Somlai, and
Vizer~\cite{KissPompeiu} used the connection with the discrete Pompeiu
problem to obtain the case $N=p^a q^2$.

For a spectral pair $(A,B)$, column orthogonality gives
\[
\sum_{x\in A}\zeta_N^{x(y-y')}=0
\qquad (y,y'\in B,\ y\ne y').
\]
Thus the additive order of $y-y'$ determines a cyclotomic divisor of
the mask polynomial of $A$. Orthogonality of the rows gives the
corresponding divisors for the mask of $B$.

For $N=p^a q^b$, Malikiosis~\cite[Sections~8--12]{Malikiosis}
reduced the problem to primitive, absorption-free spectral pairs at
a minimal counterexample modulus. He compared zeros at orders
$p^{j+1}$ for one member of the pair with zeros at orders $N/p^j$
for the other. His root deficits measure exceptions to this pairing.
He combined estimates on these deficits with divisibility constraints
on $|A|=|B|$. For $p<q$, these estimates prove the spectral-to-tiling
implication whenever $a\le9$, $b\le6$, or
$p^{a-2}<q^4$~\cite[Theorem~1.4]{Malikiosis}.

\Needspace{8\baselineskip}
We use the same reduction to a minimal counterexample. As $v$
increases, we track the proportion of points in each residue class
modulo $q$ within each nonempty class of $A$ modulo $p^v$. Estimates for the weights of connected components
in the graphs of character rows control the denominators of these
proportions. The mask zeros constrain the differences of the
proportions, and together the two restrictions yield the theorem
for arbitrary positive exponents.
\nopagebreak[4]

\begin{theorem}\label{thm:main}
	Let $p,q$ be distinct primes and let $a,b\geq 1$. Every nonempty spectral subset
	of $\Z_{p^a q^b}$ tiles $\Z_{p^a q^b}$.
\end{theorem}

Combining Theorem~\ref{thm:main} with the two-prime tiling-to-spectral
theorem gives Fuglede's conjecture on $\Z_{p^a q^b}$. For every nonempty
$A\subseteq\Z_{p^a q^b}$,
\[
A\text{ is spectral}
\quad\Longleftrightarrow\quad
A\text{ tiles }\Z_{p^a q^b}.
\]

\Needspace{10\baselineskip}
\subsection*{Malikiosis's reduction and our argument}

For $D\subseteq\Z_N$, write $M_D(X)=\sum_{x\in D}X^x$ using the
representatives $0,\ldots,N-1$ for the exponents, and put
$D[d]=M_D(\zeta_d)$ for $d\mid N$.

We continue from the minimal-counterexample reduction in
Malikiosis~\cite[Theorem~10.9]{Malikiosis}. Take a minimal modulus
$N=p^a q^b$ with $p<q$ for which the spectral-to-tiling implication
fails, and choose a nontiling spectral set $A$ of maximum cardinality
with spectrum $B$. Malikiosis shows that $B$ is itself a nontiling
spectral set of maximum cardinality. Both sets are primitive and
absorption-free.
At a $p$-level $j<a$, this means that some class modulo $p^j$ has
at least two nonempty refinements modulo $p^{j+1}$. He also obtains
$A[p]=A[q]=B[N]=B[q^b]=0$ and rules out a spectrum made up exclusively
of $p$-cycles or exclusively of $q$-cycles. A prime-order cycle is a
coset of the corresponding prime-order subgroup.

The zero $B[N]=0$ gives a decomposition into $p$- and $q$-cycles
by the theorem of Lam and Leung~\cite{LamLeung}. Malikiosis
relates zeros of the two masks by following how residue classes
split at successive prime powers. In the notation of
Section~\ref{sec:obstruction}, let $m$ be the largest $0\le j\le a$
such that $B[N/p^j]\ne0$, and put $K=m+1$. Applying his
Lemma~10.3 and Proposition~10.4 to $(B,A)$ gives
\[
B[N/p^j]=0\quad\Longleftrightarrow\quad A[p^{j+1}]=0
\quad(0\le j\le m),\qquad A[p^Kq]=0.
\]
Taking $j=m$ gives $A[p^K]\ne0$. Let $k$ be the largest index for
which $A[p^j]=0$ for all $1\le j\le k$. Then $1\le k<K\le a$.

We use these same structural conclusions to study the proportions
of points in each residue modulo $q$ within classes modulo powers of
$p$.

For a nonempty class $C$ of $A$ modulo $p^v$ and the class $A_\alpha$ with
residue $\alpha$ modulo $q$, put $\pi_\alpha(C)=|C\cap A_\alpha|/|C|$.
For $0\le v\le k$, our proof gives
\[
\pi_\alpha(C)\in\Zp,\qquad
\pi_\alpha(C)-\pi_\beta(C)\in\Z[1/p].
\]
The reduced denominator is coprime to $p$ in $\Zp$ and a power of $p$
in $\Z[1/p]$.

The root relation gives $B[N/p^j]=0$ for $0\le j<k$.
Section~\ref{sec:spectral} uses these zeros to show that the set left
after removing the complete $q$-cycles from $B$ is invariant under
translation by $N/p^k$.
This periodicity and row orthogonality give the commutation of
the relevant column-space projections with the coordinate projections
onto the row classes. Fourier transformation within each complete cycle
expresses the row compressions as weighted character projections on
$\Z_{q^{b-1}}$. The component-weight bounds in Sections~\ref{sec:rows}
and~\ref{sec:rank} allow at most one extra factor of $p$ in a denominator
at each refinement of a row class. For $X=C\cap A_\alpha$, the total weight
$w_X(\Z_{q^{b-1}})$ therefore belongs to $p^{-v}\Zp$. This weight equals
$q|X|/|A|$, and $|C|=|A|/p^v$, so multiplying it by $p^v/q$ gives
the proportion $\pi_\alpha(C)$ in $\Zp$.

For the difference of two proportions, Section~\ref{sec:residues}
subdivides $C$ into classes at successive powers of $p$. Each split
has either $p$ children of equal size or a single nonempty child.
A class contained in one residue modulo $q$ is kept as a leaf.
At depth $K$, the unique nonempty child of a mixed parent is uniform
modulo $q$ by the zero $A[p^Kq]=0$. The resulting leaves have relative
sizes that are powers of $1/p$. The difference of two proportions
vanishes on a uniform leaf and is an integer on a leaf contained in
one residue class. Taking the weighted sum over the leaves gives
the condition in $\Z[1/p]$.

The two conditions make every difference of proportions an integer.
Consequently each class through depth $k$ is uniform modulo $q$ or is
contained in one residue. The root $A$ is uniform, and $p<q$ forces
uniformity to pass to its equal-sized children at each of these depths.
At depth $k$, uniformity and $A[p^{k+1}]\ne0$ force a single nonempty
child in each class. This contradicts absorption-freeness, which
requires at least two nonempty children in some class at that depth.

\section{Fourier notation and the minimal counterexample}\label{sec:obstruction}

\subsection{Mask polynomials and coefficient comparison}

Write $\Z_n=\Z/n\Z$ and $\zeta_n=\exp(2\pi i/n)$. For a set or indexed
list $D$ in $\Z_n$, define
\[
M_D(X)=\sum_{z\in D}X^z,\qquad D[d]=M_D(\zeta_d)\quad(d\mid n).
\]
The sum uses representatives in $\{0,\ldots,n-1\}$ and counts each list
index once. Since $\zeta_d^n=1$, the value $D[d]$ is independent of the
representatives. If $D[d]=0$, then $M_D$ vanishes at every primitive
$d$th root, since these roots are conjugate over $\Q$.

We use the Fourier convention
\[
\widehat f(y)=\sum_{x\in\Z_n}f(x)\zeta_n^{xy},
\qquad
f(x)=\frac1n\sum_{y\in\Z_n}\widehat f(y)\zeta_n^{-xy}.
\]
Thus $\widehat{\one_D}(y)=M_D(\zeta_n^y)$ for a set $D$.
The order of $\zeta_n^y$ is the additive order
$\ord(y)=n/\gcd(n,y)$. For a spectral pair $(A,B)$ in $\Z_N$,
orthogonality of the columns and rows of
$|A|^{-1/2}(\zeta_N^{xy})_{x\in A,\,y\in B}$ gives
\begin{equation}\label{eq:spectral}
	\begin{aligned}
		y,y'\in B,\quad y\ne y'&\quad\Longrightarrow\quad A[\ord(y-y')]=0,\\
		x,x'\in A,\quad x\ne x'&\quad\Longrightarrow\quad B[\ord(x-x')]=0.
	\end{aligned}
\end{equation}
The transposed matrix is also unitary, so $(B,A)$ is spectral as well.

\Needspace{6\baselineskip}
\begin{lemma}\label{thm:coeff}
	Let $\ell$ be a prime, $s\ge1$, and $D_0=\ell^{s-1}$. Let $E$ be a subfield
	of $\mathbb C$ over which $\Phi_{\ell^s}$ is irreducible, and let $\xi$ be a
	primitive $\ell^s$th root of unity. For coefficients $c_{r,j}\in E$,
	\begin{equation}\label{eq:coeff}
		\sum_{r=0}^{D_0-1}\sum_{j=0}^{\ell-1}c_{r,j}\xi^{\,r+jD_0}=0
		\quad\Longleftrightarrow\quad
		c_{r,0}=\cdots=c_{r,\ell-1}\quad(0\le r<D_0).
	\end{equation}
\end{lemma}
\begin{proof}
	Let $P(X)$ be the polynomial on the left of \eqref{eq:coeff}, with $X$
	in place of $\xi$. Since $\Phi_{\ell^s}$ is the minimal polynomial of
	$\xi$ over $E$, the equality $P(\xi)=0$ is equivalent to
	$\Phi_{\ell^s}\mid P$ in $E[X]$. We have
	\[
	\Phi_{\ell^s}(X)=\sum_{j=0}^{\ell-1}X^{jD_0},
	\qquad
	\deg P\le \ell^s-1.
	\]
	The quotient therefore has degree at most $D_0-1$ and can be written as
	$\sum_{r=0}^{D_0-1}b_rX^r$. Its product with $\Phi_{\ell^s}$ has
	coefficient $b_r$ at every position $r+jD_0$. This gives the equalities in
	\eqref{eq:coeff}. Conversely, when those equalities hold, taking
	$b_r=c_{r,0}$ gives that factorization of $P$ and hence $P(\xi)=0$.
\end{proof}

For a mask in a group of order divisible by $\ell^s$, the coefficients
in \eqref{eq:coeff} are obtained by first reducing the exponents modulo
$\ell^s$. In particular, $c_{r,j}$ counts the indices in the residue
$r+jD_0$, so vanishing is equivalent to equal counts at the $\ell$
positions with common residue $r$ modulo $D_0$.

For coprime positive integers $u,v$, B\'ezout's identity gives
$\Q(\zeta_u,\zeta_v)=\Q(\zeta_{uv})$. The cyclotomic degree formula gives
\[
[\Q(\zeta_u,\zeta_v):\Q(\zeta_v)]
=\frac{\varphi(uv)}{\varphi(v)}=\varphi(u).
\]
Thus $\Phi_u$ is irreducible over $\Q(\zeta_v)$. For distinct primes
$p,q$, this lets us use lemma~\ref{thm:coeff} at $p$-power roots over
$\Q(\zeta_q)$ and at primitive $p$th roots over $\Q(\zeta_{q^e})$
for every $e\ge0$, in addition to its use over $\Q$.

\Needspace{5\baselineskip}
\subsection{The minimal counterexample and its mask zeros}
Suppose the spectral-to-tiling implication fails for cyclic groups of
order $N_0=p^{a_0}q^{b_0}$ with $p<q$ and $a_0,b_0\ge1$. Among the divisors of
$N_0$ for which the implication fails, choose $N=p^a q^b$ minimal under
divisibility. The implication holds for the trivial group and for
cyclic groups of prime-power order, so $a,b\ge1$. Choose a nontiling spectral set $A$ of maximum
cardinality in $\Z_N$, and let $B$ be a spectrum of $A$.

Malikiosis~\cite[Theorem~10.9]{Malikiosis} gives the following properties.
\begin{enumerate}[label=\textup{(\roman*)},leftmargin=*]
	\item $B$ is also a nontiling spectral set of maximum cardinality, and
	$A[p]=A[q]=B[N]=B[q^b]=0$.
	\item $B$ cannot be written exclusively as a union of $p$-cycles or exclusively as
	a union of $q$-cycles. An $\ell$-cycle is a coset of the subgroup of
	order $\ell$.
	\item For every $0\le j<a$, at least one class of $A$ modulo $p^j$ has two or more nonempty children  modulo $p^{j+1}$. This is the absorption-free property
	at the $p$-levels.
\end{enumerate}

Set
\[
m=\max\{j\in\{0,\ldots,a\}:B[N/p^j]\ne0\}.
\]
If this index set were empty, $\widehat{\one_B}(y)=0$ for every $y$ not
divisible by $q$, since $\ord(y)$ is then one of the orders $N/p^j$.
For $q\mid y$, the character is invariant under the shift $N/q$.
Fourier inversion would give $\one_B(x+N/q)=\one_B(x)$, making $B$ a
union of $q$-cycles, contrary to (ii). Thus $m$ exists, and $m<a$
because $B[q^b]=0$.

The pair $(B,A)$ is spectral and $B$ is nontiling at the minimal
modulus $N$. Thus Malikiosis's Lemma~10.3 and Proposition~10.4
\cite{Malikiosis}, applied with $B$ first, give
\begin{equation}\label{eq:pattern}
	B[N/p^j]=0\quad\Longleftrightarrow\quad A[p^{j+1}]=0
	\quad(0\le j\le m),\qquad A[p^{m+1}q]=0 .
\end{equation}
Put $K=m+1$. At $j=m$, the definition of $m$ gives $A[p^K]\ne0$ by
\eqref{eq:pattern}. and let $k$ be the number of consecutive
zeros in $A[p],A[p^2],\ldots$ starting at $p$. Since $A[p]=0$ and $A[p^K]\ne0$,
we have $1\le k<K\le a$. Using \eqref{eq:pattern} for $0\le j<k$
then yields the zeros of $B$ in
\begin{equation}\label{eq:firstzeros}
	\begin{gathered}
		1\le k<K\le a,\qquad
		A[p^j]=0\ (1\le j\le k),\qquad A[p^{k+1}]\ne0,\\
		B[N/p^j]=0\ (0\le j<k),\qquad
		A[p^K]\ne0,\qquad A[p^Kq]=0.
	\end{gathered}
\end{equation}
% END SECTION 1

\section{Weighted character rows}\label{sec:rows}

\subsection{The graph and its weights}

In Sections~\ref{sec:rows} and \ref{sec:rank}, $n$ is a power of $q$, including
$1$. Let $F=(z_\xi)_{\xi\in I}$ be an indexed list in $\Z_n$, with $f=|I|>0$.
Repeated frequency values have distinct indices. For $x\in\Z_n$, let
$u_x=(\zeta_n^{xz_\xi})_{\xi\in I}$ be a row vector. Let $w:\Z_n\to\Q_{\ge0}$. Throughout this discussion we assume that
\begin{equation}\label{eq:gram}
	P=P(F,w)=\sum_{x\in\Z_n}w(x)u_x^*u_x
\end{equation}
is an orthogonal projection, so $P=P^*=P^2$. The notation $u_x^*$ denotes the conjugate
transpose, so $u_x^*u_x$ is an $f\times f$ matrix. For $C\subseteq\Z_n$ write
$w(C)=\sum_{x\in C}w(x)$, and put $\sigma=w(\Z_n)$.

Each entry of $P$ and each inner product of rows can be written as
\begin{equation}\label{eq:gramsentries}
 P_{\xi\eta}=\sum_x w(x)\zeta_n^{x(z_\eta-z_\xi)},\qquad
 u_xu_y^*=\sum_{\xi\in I}\zeta_n^{(x-y)z_\xi}
 =M_F(\zeta_n^{x-y}),
\end{equation}
where $M_F(X)=\sum_{\xi\in I}X^{z_\xi}$ is the mask of the indexed
list. In particular $P_{\xi\xi}=\sigma$ and the entries belong to $\Q(\zeta_n)$.
Let $\supp w=\{x\in\Z_n:w(x)>0\}$ and set
$
V=\operatorname{span}\{u_x:x\in\supp w\}.
$
For any row $v$,
\[
 vPv^*=\sum_x w(x)|vu_x^*|^2.
\]
As $P=P^*=P^2$, the left side equals $\|vP\|^2$. It is zero exactly
when $v$ is orthogonal to all $u_x$ with $w(x)>0$, so $\ker P=V^\perp$. It follows that
$
\operatorname{Ran}P=(\ker P)^\perp=(V^\perp)^\perp=V,
$
where the last equality uses finite dimensionality. The restriction
of orthogonal projection $P$ to its range $V$ is the identity operator. In particular, if
$x\in\supp w$, then $u_x\in V$ and hence $u_xP=u_x$.

Form a graph with vertex set $\supp w$ by joining distinct vertices $x$ and $y$
when $u_xu_y^*\ne0$. Thus the rows in different connected components are orthogonal,
since a nonzero inner product would be an edge joining the components.
The components partition $\supp w$, and the weights outside
$\supp w$ are zero. It follows that
\begin{equation}\label{eq:component-weights}
	\sigma=w(\Z_n)=\sum_{C} w(C),
\end{equation}
where the sum is over all connected components of the graph.

Write
\[
\Zp=\{a/b:a,b\in\Z,\ b\ne0,\ p\nmid b\}.
\]
These are the rationals whose reduced denominators are prime to $p$.
For an integer $d\ge0$, say that $w$ satisfies $\Rp_d$ on $F$ if
\begin{equation}\label{eq:mass}
	w(C)\in p^{-d}\Zp\qquad\text{for every connected component }C .
\end{equation} 
Note that $p^{-d}\Zp$ is closed under addition. Thus if $w$ satisfies $\Rp_d$, equations
\eqref{eq:component-weights} and \eqref{eq:mass} give
$\sigma\in p^{-d}\Zp$.  When $P=0$, the diagonal identity
$P_{\xi\xi}=\sigma$ and nonnegativity of the weights give $w=0$.
The graph is then empty, so $\Rp_d$ holds vacuously.

\Needspace{6\baselineskip}
\begin{lemma}\label{lem:orthogonal-sums}
	Let $P_i=P(F,w_i)$ for $1\le i\le s$ be orthogonal projections
	on the same indexed list $F$, with nonnegative rational weights $w_i$.
	Assume $P_iP_j=0$ for $i\ne j$.
	Set $w=\sum_i w_i$. Then $P(F,w)=\sum_i P_i$ is a projection,
	and its graph is the disjoint union of the graphs for the weights
	$w_i$. If $C$ is a component from $w_i$, then $w(C)=w_i(C)$.
	In particular, $w$ satisfies $\Rp_d$ if and only if every $w_i$
	satisfies $\Rp_d$.
\end{lemma}
\begin{proof}
	The expression \eqref{eq:gram} is linear in $w$, so
	$P(F,w)=\sum_iP_i$. This sum is self-adjoint, and
	\[
	\left(\sum_iP_i\right)^2
	=\sum_iP_i^2+\sum_{i\ne j}P_iP_j
	=\sum_iP_i .
	\]
	
	Take $x\in\supp w_i$ and $y\in\supp w_j$ with $i\ne j$. By the
	row-space description, $u_xP_i=u_x$ and $u_yP_j=u_y$. The latter
	identity also gives $P_j u_y^*=u_y^*$ by taking adjoints. Thus
	\[
	u_xu_y^*=(u_xP_i)(P_j u_y^*)=u_xP_iP_j u_y^*=0.
	\]
	If the supports met at a point $x$, taking $y=x$ would give
	$u_xu_x^*=0$, contrary to $u_xu_x^*=f>0$. Thus the supports
	are disjoint. Nonnegativity of the weights gives
	$\supp w=\bigcup_i\supp w_i$. Let $G$ and $G_i$ be the graphs of $w$ and $w_i$ on the common list $F$.
	Since $u_xu_y^*=0$ whenever $x\in\supp w_i$, $y\in\supp w_j$, and $i\ne j$, which excludes all
	edges between distinct supports. Thus
	\[
	E(G)
	=\bigsqcup_{i=1}^{s}
	\bigl\{\{x,y\}:x,y\in\supp w_i,\ x\ne y,\ u_xu_y^*\ne0\bigr\}
	=\bigsqcup_{i=1}^{s}E(G_i).
	\]
	Together with $\supp w=\bigsqcup_{i=1}^{s}\supp w_i$, this gives
	$G=\bigsqcup_{i=1}^{s}G_i$.
	
	On a component $C$ from the support of $w_i$, every other weight
	function vanishes, so $w(x)=w_i(x)$ for $x\in C$ and
	$w(C)=w_i(C)$. Since every component of the sum comes from exactly
	one summand, \eqref{eq:mass} for $w$ is equivalent to
	\eqref{eq:mass} for all the $w_i$.
\end{proof}
% END SECTION 2.1

\subsection{Changing the ambient group}

Let $L\mid n$ and suppose $\supp w$ lies in  $x_0+(n/L)\Z_n$. The map $\varphi(t)=x_0+(n/L)t$ is a bijection from $\Z_L$ onto this residue class.
Set $\widetilde w(t)=w(\varphi(t))$ and
$\widetilde F=(\bar z_\xi)_{\xi\in I}$, where $\bar z_\xi=z_\xi\bmod L$
and each index $\xi$ is retained. Then $\widetilde F[r]=F[r]$ for
$r\mid L$, since $\zeta_r^L=1$.

Let $v_t=(\zeta_L^{t\bar z_\xi})_{\xi\in I}$ be the character row for
$\widetilde F$ and set $D_{x_0}=\diag_{\xi\in I}(\zeta_n^{x_0z_\xi})$.
We obtain
\begin{equation}\label{eq:fold}
	u_{\varphi(t)}
	=\bigl(\zeta_n^{(x_0+(n/L)t)z_\xi}\bigr)_{\xi\in I}
	=\bigl(\zeta_L^{t\bar z_\xi}\zeta_n^{x_0z_\xi}\bigr)_{\xi\in I}
	=\bigl(\zeta_L^{t\bar z_\xi}\bigr)_{\xi\in I}
	\diag_{\xi\in I}\bigl(\zeta_n^{x_0z_\xi}\bigr)
	=v_tD_{x_0}.
\end{equation}
The matrix $D_{x_0}$ is unitary. Summing over the residue class gives
\[
P(F,w)=\sum_{t\in\Z_L}\widetilde w(t)(v_tD_{x_0})^*(v_tD_{x_0})
=D_{x_0}^*P(\widetilde F,\widetilde w)D_{x_0}.
\]
Hence $P(\widetilde F,\widetilde w)$ is an orthogonal projection.

The map $\varphi$ takes $\supp\widetilde w$ bijectively onto $\supp w$.
For $t,s\in\supp\widetilde w$,
\[
u_{\varphi(t)}u_{\varphi(s)}^*
=v_t D_{x_0}D_{x_0}^*v_s^*
=v_t v_s^* .
\]
Thus $\varphi$ preserves adjacency in both directions, it induces
a one-to-one correspondence between the connected components of the two graphs,
$
\widetilde C \longleftrightarrow C=\varphi(\widetilde C).
$
For each such pair,
\[
\widetilde w(\widetilde C)
=\sum_{t\in\widetilde C}w(\varphi(t))
=\sum_{x\in C}w(x)
=w(C).
\]
In particular, $\widetilde w(\Z_L)=w(\Z_n)$, and for every $d\ge0$,
\[
w\text{ satisfies }\Rp_d\text{ on }F
\quad\Longleftrightarrow\quad
\widetilde w\text{ satisfies }\Rp_d\text{ on }\widetilde F .
\]

\Needspace{6\baselineskip}
\begin{lemma}\label{lem:coset}
	Let $Y$ be a nonempty subset of $\Z_n$ and let
	$L=\max_{x,y\in Y}\ord(x-y)$, where $\ord$ denotes additive order.
	Then $L$ is the least cardinality of a coset containing $Y$.
	If $Y=\bigcup_{i=1}^s Y_i$ with $s\ge2$ and each $Y_i$ nonempty,
	then there exist $i\ne j$ and points $x\in Y_i$, $y\in Y_j$ such that
	$\ord(x-y)=L$. In particular, if the row graph of $P(F,w)$ has at least two connected
	components and $L$ is the least cardinality of a coset containing
	$\supp w$, then $F[L]=0$.
\end{lemma}
\begin{proof}
	For $r\mid n$ and $h\in\Z_n$,
	\[
	\ord(h)\mid r
	\quad\Longleftrightarrow\quad rh=0\text{ in }\Z_n
	\quad\Longleftrightarrow\quad h\equiv0\pmod{n/r}.
	\]
	Every order in $\Z_n$ is a power of $q$, including $\ord(0)=1$.
	Fix $x_0\in Y$. Since $\ord(x-x_0)\le L$, it divides $L$, so
	$x\equiv x_0\pmod{n/L}$ for every $x\in Y$. This class has size $L$.
	Conversely, choose $x_*,y_*\in Y$ with $\ord(x_*-y_*)=L$.
	If $Y$ lies in a coset of size $r$, then $x_*$ and $y_*$ lie in the
	same residue class modulo $n/r$. Hence
	\[
	x_*\equiv y_*\pmod{n/r}
	\quad\Longrightarrow\quad
	r(x_*-y_*)=0 \text{ in }\Z_n
	\quad\Longrightarrow\quad
	L=\ord(x_*-y_*)\mid r.
	\]
	Thus every coset containing $Y$ has size at least $L$.
	
	Suppose $\ord(x-y)<L$ whenever $x\in Y_i$, $y\in Y_j$, and
	$i\ne j$. Each such order divides $L/q$, so $x\equiv y\pmod{qn/L}$
	for all these pairs.
	For $x,x'\in Y_i$, choose $y\in Y_j$ with $j\ne i$. Then
	$
	x\equiv y\equiv x'\pmod{qn/L}.
	$
	Thus all of $Y$ lies in a residue class of size $L/q$, contradicting
	the first assertion.
	
	For the last assertion, take the connected components of $\supp w$
	as the sets $Y_i$. The second assertion gives $x,y$ in different components with $\ord(x-y)=L$.
	There is no edge between them, hence
	\[
	0=u_xu_y^*=M_F(\zeta_n^{x-y}).
	\]
	The root $\zeta_n^{x-y}$ has order $\ord(x-y)=L$. Since $M_F$ has
	integer coefficients and the primitive $L$th roots are conjugate
	over $\Q$, it follows that $F[L]=M_F(\zeta_L)=0$.
\end{proof}
% END SECTION 2.2

\subsection{The \texorpdfstring{$q$}{q}-point Fourier splitting}

\begin{lemma}\label{lem:split}
Suppose $L\mid n$, $L>1$, $F[L]=0$, and
$\supp w\subseteq x_0+(n/L)\Z_n$. By Section~2.2, work in $\Z_L$ and
retain the notation $F,w,u_x$ and $P=P(F,w)$. Put $T=L/q$.
The columns can be matched into groups with frequencies $z+jT$
for $0\le j<q$ and $0\le z<T$. Let $F'$ in $\Z_T$ record one base
frequency $z$ per group. For any such matching, $P$ is unitarily
equivalent to the direct sum of orthogonal projections
\[
\bigoplus_{s=0}^{q-1}Q_s,\qquad Q_s=P(F',w_s),\qquad
w_s(t)=q\,w(s+qt)\quad(t\in\Z_T).
\]
The row components correspond bijectively to the components of the
channel graphs, with $w_s(C')=q\,w(C)$ for a component $C$ and its
image $C'$ in channel $s$. Thus $w$ satisfies $\Rp_d$ on $F$ if and
only if every $w_s$ satisfies it on $F'$.
\end{lemma}

\begin{proof}
\textit{Matching the columns.}
Since $F[L]=0$, Lemma~\ref{thm:coeff} shows that for each $0\le z<T$,
the frequencies $z+jT$ with $0\le j<q$ have the same multiplicity $c_z$.
Match these columns into $c_z$ groups so that the $q$ columns
in each group correspond one-to-one to the $q$ frequencies.

Let $\mathcal G$ be the set of matched groups. For $g\in\mathcal G$,
write $z_g$ for its base frequency and $\xi(g,j)\in I$ be the original index of the
column in group $g$ with frequency $z_g+jT$.. Then
\[
z_{\xi(g,j)}=z_g+jT,\qquad
F'=(z_g)_{g\in\mathcal G},\qquad
|\mathcal G|=\sum_{z=0}^{T-1}c_z=f/q .
\]

\smallskip
\noindent\textit{The column change.}
Write each $x\in\Z_L$ uniquely as $x=s+qt$ with $0\le s<q$ and
$t\in\Z_T$. We construct $U=RHSD$, acting by right multiplication
on each character row.

\emph{Group order.}
The permutation $R$ puts the columns in $(g,j)$ order, with
$(aR)_{g,j}=a_{\xi(g,j)}$ for a row $a$. Hence
\begin{equation}\label{eq:rowDFT}
(u_{s+qt}R)_{g,j}
=\zeta_L^{(s+qt)(z_g+jT)}
=\zeta_L^{sz_g}\zeta_T^{tz_g}\zeta_q^{sj}.
\end{equation}
Here $\zeta_L^q=\zeta_T$, $\zeta_L^T=\zeta_q$, and
$qt\,jT=tjL$ contributes a factor of $1$.

\emph{Fourier transform within each group.}
Put $H=\bigoplus_{g\in\mathcal G} H_q$, where
$(H_q)_{j,k}=q^{-1/2}\zeta_q^{-jk}$ for $0\le j,k<q$.
The geometric sum gives
$(H_q^*H_q)_{k,k'}=q^{-1}\sum_{j=0}^{q-1}\zeta_q^{j(k-k')}
=\delta_{k,k'}$, so $H_q$ and $H$ are unitary. It also gives
\[
\begin{aligned}
(u_{s+qt}RH)_{g,k}
 &=\frac{\zeta_L^{sz_g}\zeta_T^{tz_g}}{\sqrt q}
   \sum_{j=0}^{q-1}\zeta_q^{(s-k)j}\\
 &=\sqrt q\,\delta_{k,s}\zeta_L^{sz_g}\zeta_T^{tz_g}.
\end{aligned}
\]

\emph{Channel order and phase correction.}
The permutation $S$ changes $(g,k)$ to $(k,g)$, with
$(bS)_{k,g}=b_{g,k}$. In this order, let
$D=\diag_{(k,g)}(\zeta_L^{-kz_g})$. Its entries on the diagonal
have modulus $1$, so $U$ is unitary, and
\[
(u_{s+qt}U)_{k,g}
=\sqrt q\,\delta_{k,s}\zeta_L^{(s-k)z_g}\zeta_T^{tz_g}
=\sqrt q\,\delta_{k,s}\zeta_T^{tz_g}.
\]
Thus the transformed row has $\sqrt q\,u'_t$ in channel $s$ and
zero in the other channels, where
$u'_t=(\zeta_T^{tz_g})_{g\in\mathcal G}$ is the row for $F'$.

\smallskip
\noindent\textit{Channel matrices.}
The row $u_{s+qt}U$ has $\sqrt q\,u'_t$ in channel $s$ and zeros
elsewhere. Thus
\[
\begin{aligned}
	U^*PU
	&=U^*\left(\sum_{x\in\Z_L}w(x)u_x^*u_x\right)U\\
	&=\sum_{s=0}^{q-1}\sum_{t\in\Z_T}
	w(s+qt)(u_{s+qt}U)^*(u_{s+qt}U)\\
	&=\bigoplus_{s=0}^{q-1}
	\left(\sum_{t\in\Z_T}q\,w(s+qt)(u'_t)^*u'_t\right).
\end{aligned}
\]
The $s$th block is $Q_s=P(F',w_s)$ with $w_s(t)=q\,w(s+qt)$. Since $U^*PU$ is an orthogonal projection, each $Q_s$ is one as well.

\smallskip
\noindent\textit{Connected components and weights.}
For $x=s+qt$ and $y=r+qh$, unitarity gives
\[
u_xu_y^*=(u_xU)(u_yU)^*=q\,\delta_{s,r}u'_t(u'_h)^* .
\]
The positive-weight points are in bijection under
\[
\bigsqcup_{s=0}^{q-1}(\{s\}\times\supp w_s)
\longrightarrow \supp w,\qquad (s,t)\longmapsto s+qt.
\]
The inner product formula preserves adjacency in each channel and
gives zero inner products between channels. Thus $\psi_s(t)=s+qt$
takes each channel component $C'$ to a component $C=\psi_s(C')$,
and all components arise in this way. Their weights satisfy
\[
w_s(C')=q\sum_{t\in C'}w(s+qt)=q\,w(C).
\]
Since $q$ is a unit in $\Zp$, this yields the equivalence of
$\Rp_d$ for $w$ and for all the $w_s$.
\end{proof}

\subsection{Repeated frequencies and fixed rows}

For an occurring frequency $z$, write $I_z=\{\xi\in I:z_\xi=z\}$,
$c_z=|I_z|$, and let $e_z=c_z^{-1/2}\one_{I_z}$ be the corresponding
unit row. Write $P[I_z,I_z]=(P_{\xi\eta})_{\xi,\eta\in I_z}$ for the principal
submatrix obtained by keeping the rows and columns indexed by $I_z$.
For $\xi,\eta\in I_z$ we have $z_\xi=z_\eta=z$, so
\eqref{eq:gramsentries} gives $P_{\xi\eta}=\sum_{x\in\Z_n}w(x)=\sigma$. 
Since $I_f-P$ is an orthogonal projection, we have
\[
\begin{aligned}
e_zPe_z^*
 &=\frac{1}{c_z}\sum_{\xi,\eta\in I_z}P_{\xi\eta}
   =c_z\sigma,\\
1-c_z\sigma
 &=e_z(I_f-P)e_z^*=\|e_z(I_f-P)\|^2.
\end{aligned}
\]
Consequently
\begin{equation}\label{eq:saturation}
e_zPe_z^*=c_z\sigma\le1,\qquad
e_zP=e_z\quad\Longleftrightarrow\quad c_z\sigma=1 .
\end{equation}
We call $z$ \emph{saturated} when $c_z\sigma=1$.

\Needspace{6\baselineskip}

\begin{lemma}\label{lem:saturated}
Suppose a nonzero projection $P=P(F,w)$ has a saturated frequency.
For every connected component $C$ of its row graph there is an integer
$r_C\ge0$ such that $w(C)=q^{-r_C}\sigma$. In particular,
$\sigma\in p^{-d}\Zp$ implies $\Rp_d$, and weights giving the identity
projection $I_f$ satisfy $\Rp_0$.
\end{lemma}

\begin{proof}
We use strong induction on the exponent of $q$ in $n$. If the graph is connected,
its only component has weight $\sigma$. When $n=1$, the nonzero projection has $\supp w=\{0\}$, so the row graph is connected. This covers $n=1$

Suppose the graph has at least two components, and fix a saturated
frequency $z_*$ of multiplicity $c$, so $c\sigma=1$.
Let $L$ be the size of the smallest coset containing the support.
Then $L>1$, and Lemma~\ref{lem:coset} gives $F[L]=0$. By Section~2.2
we work in $\Z_L$ with the same notation for the model and the same
total and component weights.

Let $\bar z_*$ be the reduction of $z_*$ modulo $L$, of multiplicity
$c'$ in the reduced list. All $c$ original column indices occur at
$\bar z_*$ after reduction. Applying \eqref{eq:saturation} in the
reduced model gives
\[
c'\ge c,\qquad c'\sigma\le1=c\sigma.
\]
Thus $c'=c$ and $\bar z_*$ is saturated.

Put $T=L/q$ and write $\bar z_*=z_0+j_0T$ with $0\le z_0<T$ and
$0\le j_0<q$. In Lemma~\ref{lem:split}, the frequencies $z_0+jT$
have equal multiplicities for $0\le j<q$. The position $j_0$ has
multiplicity $c$, so the $q$ positions form $c$ matched groups.
By Lemma~\ref{lem:split}, every channel uses the same short list $F'$.
The $c$ groups of base frequency $z_0$ contribute one copy each to
$F'$, so $z_0$ has multiplicity $c$ in every channel.

Let $\sigma_s=\sum_{t\in\Z_T}w_s(t)$ be the total weight in channel
$s$. By \eqref{eq:saturation} and $w_s(t)=q\,w(s+qt)$,
\[
c\sigma_s\le1\quad(0\le s<q),\qquad
\sum_{s=0}^{q-1}c\sigma_s
=cq\sum_{s=0}^{q-1}\sum_{t\in\Z_T}w(s+qt)
=cq\sigma=q .
\]
Each of the $q$ terms must therefore equal $1$. Thus
$c\sigma_s=1$ and $\sigma_s=\sigma>0$ for all $s$, so every $Q_s$
is nonzero with $z_0$ saturated.
As $T=L/q<n$, induction applies to each of these projections.

For any original component $C$, let $C'$ be the corresponding
component in channel $s$. By induction and Lemma~\ref{lem:split},
there is an integer $r\ge0$ with
\[
q\,w(C)=w_s(C')=q^{-r}\sigma_s=q^{-r}\sigma ,
\qquad w(C)=q^{-(r+1)}\sigma .
\]
The $\Rp_d$ assertion follows since $q$ is a unit in $\Zp$.
For $P=I_f$, \eqref{eq:gramsentries} gives $\sigma=1$, and
\eqref{eq:saturation} gives $c_z=1$ at every occurring frequency.
Hence each such frequency is saturated and the result applies with
$d=0$.
\end{proof}

\Needspace{6\baselineskip}
\begin{lemma}\label{lem:digits}
Let $P=P(F,w)$ be on $\Z_T$ with $T=q^k$, and take the frequency
representatives $0\le z_\xi<T$. Suppose $h>1$ is a power of $q$
and $a_0\in\Z$ is prime to $q$. If a row $v$ satisfies $vP=v$ and
\begin{equation}\label{eq:digitform}
	v_\xi=\zeta_{hT}^{a_0z_\xi}b_{z_\xi},
	\qquad b_{z_\xi}\in\Q(\zeta_h)
	\quad(\xi\in I),
\end{equation}
then every occurring frequency $z$ with $b_z\ne0$ is saturated.
\end{lemma}

\begin{proof}
Let $r$ be the extension by zero of a row $r_J$ on a set $J\subseteq I$.
Since $(rP)|_J=r_JP[J,J]$, the equation $rP=r$ implies
$r_JP[J,J]=r_J$. Conversely, if $r_JP[J,J]=r_J$, the projection
identity for $I_f-P$ gives
\[
\|r(I_f-P)\|^2
=r(I_f-P)r^*
=r_J\bigl(I_{|J|}-P[J,J]\bigr)r_J^*=0,
\]
so $rP=r$.

For $0\le j\le k$ and $0\le c<q^j$, set
\[
J_j(c)=\{\xi\in I:z_\xi\equiv c\pmod{q^j}\}.
\]
Let $v^{[j,c]}$ be equal to $v$ on $J_j(c)$ and zero outside it.
We show by induction that $v^{[j,c]}P=v^{[j,c]}$ for every $c$.
At $j=0$, this is $vP=v$. Suppose $j<k$ and the assertion holds.
Fix a nonempty class $J=J_j(c)$. The observation above gives
\[
v_JP[J,J]=v_J,\qquad v_J=v|_J .
\]

Put $s=T/q^j$. For $\xi\in J$, write
$z_\xi=c+q^j t_\xi$ with $0\le t_\xi<s$. For $0\le e<q$, let
\[
J_e=\{\xi\in J:a_0t_\xi\equiv e\pmod q\}.
\]
Since $a_0$ is invertible modulo $q$, there is a unique integer
$0\le t_e<q$ with $a_0t_e\equiv e\pmod q$. For $\xi\in J$, we have
$z_\xi-(c+q^j t_e)=q^j(t_\xi-t_e)$, so
\[
a_0t_\xi\equiv e\pmod q
\quad\Longleftrightarrow\quad
t_\xi\equiv t_e\pmod q
\quad\Longleftrightarrow\quad
z_\xi\equiv c+q^j t_e\pmod{q^{j+1}}.
\]
Any index satisfying the last congruence also satisfies
$z_\xi\equiv c\pmod{q^j}$ and hence belongs to $J$. Therefore
$J_e=J_{j+1}(c+q^j t_e)$.

Set $A_c=\zeta_{hT}^{a_0c}$ and $\eta=\zeta_{hs}$. For $\xi\in J$, substitute $z_\xi=c+q^jt_\xi$ into
\eqref{eq:digitform}. Dividing by $A_c$ gives
\[
\begin{aligned}
	A_c^{-1}v_\xi
	=\zeta_{hT}^{-a_0c}\zeta_{hT}^{a_0(c+q^jt_\xi)}b_{z_\xi}
	=(\zeta_{hT}^{q^j})^{a_0t_\xi}b_{z_\xi}
	=\eta^{a_0t_\xi}b_{z_\xi},
\end{aligned}
\]
and the row $A_c^{-1}v_J$ is fixed by $P[J,J]$.
For $\xi,\kappa\in J$, \eqref{eq:gramsentries} gives
\[
P_{\xi\kappa}
=\sum_{x\in\Z_T}w(x)\zeta_T^{xq^j(t_\kappa-t_\xi)}
=\sum_{x\in\Z_T}w(x)\zeta_s^{x(t_\kappa-t_\xi)}
\in\Q(\zeta_s).
\]
Put $K_0=\Q(\eta^q)$. Both $h$ and $s$ are divisible by $q$, so
\[
(\eta^q)^{s/q}=\zeta_h,\qquad
(\eta^q)^{h/q}=\zeta_s.
\]
Thus the entries of $P[J,J]$ and the coefficients $b_z$ lie in $K_0$.
Write $hs=q^\ell$ with $\ell\ge2$, then $\eta=\zeta_{q^\ell}$ and
$K_0=\Q(\eta^q)=\Q(\zeta_{q^{\ell-1}})$. The cyclotomic degree formula gives
\[
\begin{aligned}
	[\Q(\eta):K_0]
	=\frac{[\Q(\eta):\Q]}{[K_0:\Q]}
	=\frac{\varphi(q^\ell)}{\varphi(q^{\ell-1})}=q
\end{aligned}
\]
Since $\Q(\eta)=K_0(\eta)$, the minimal polynomial of $\eta$ over
$K_0$ has degree $q$. Thus $1,\eta,\ldots,\eta^{q-1}$ is a basis of $\Q(\eta)$ over $K_0$.

For each $e$, write $r^{(e)}$ for the row equal to $v$ on $J_e$ and zero elsewhere. Let $a^{(e)}=(A_c\eta^e)^{-1}(r^{(e)}|_J)$.
For $\xi\in J_e$, the exponent $(a_0t_\xi-e)/q$ is an integer, so
$a^{(e)}_\xi=(\eta^q)^{(a_0t_\xi-e)/q}b_{z_\xi}\in K_0$. The other
coordinates are zero. The sets $J_e$ partition $J$, so
\[
A_c^{-1}v_J=\sum_{e=0}^{q-1}\eta^e a^{(e)}.
\]
Substituting this into the fixed-row equation for $A_c^{-1}v_J$
gives
\[
\sum_{e=0}^{q-1}\eta^e\bigl(a^{(e)}P[J,J]-a^{(e)}\bigr)=0.
\]
Every coordinate in parentheses belongs to $K_0$. Comparing
coefficients in the basis yields $a^{(e)}P[J,J]=a^{(e)}$ for each
$e$. Multiplying by $A_c\eta^e$ gives the fixed-row equation for
$r^{(e)}|_J$. The zero-extension observation gives $r^{(e)}P=r^{(e)}$.

As $e$ varies, $t_e$ runs over $0,\ldots,q-1$, so this treats all
the classes modulo $q^{j+1}$ within $J_j(c)$. Letting $c$ vary
treats every class at level $j+1$, including empty classes, which
give zero rows. This completes the induction.

At $j=k$, the class $J_k(z)$ is exactly $I_z$. If $z$ occurs and
$b_z\ne0$, the row
\[
v^{[k,z]}=\sqrt{c_z}\,\zeta_{hT}^{a_0z}b_z\,e_z
\]
is a nonzero multiple of $e_z$ and is fixed by $P$. Hence $e_zP=e_z$,
so $z$ is saturated by \eqref{eq:saturation}. For $T=1$ the initial
level is already the final level.
\end{proof}

\section{Weight divisibility for partitioned columns}\label{sec:rank}

%The application compares a row class with its $p$ children.
%We need divisibility for individual components because the proof
%passes to individual Fourier channels.
%The condition $\Rp_d$ keeps track of the weight of each component,
%so Lemma~\ref{lem:split} allows passage to an individual channel.

In the application, we split a row class into $p$ subclasses.
We work with one Fourier channel at a time.
By Lemma~\ref{lem:split}, a component of weight $w(C)$ corresponds
to one of weight $q\,w(C)$ in its channel. Since $q$ is a unit in
$\Zp$, the condition $\Rp_d$ passes to each channel.

\subsection{A block-diagonal projection}

\begin{lemma}\label{lem:parts}
Let $P=P(F,w)$ be an orthogonal projection. Suppose it is block
diagonal for a partition $I=\bigsqcup_\beta I_\beta$ of its column
indices. If $w$ satisfies $\Rp_d$ on the sublist for a nonempty
part $I_{\beta_0}$, then $w$ satisfies $\Rp_d$ on $F$.
\end{lemma}

\begin{proof}
We use strong induction on the exponent of $q$ in $n$.
For each nonempty part $I_\beta$, write $F_\beta=(z_\xi)_{\xi\in I_\beta}$.
Restrict each character row to $I_\beta$. Then
\[
P[I_\beta,I_\beta]
=\sum_{x\in\Z_n}w(x)(u_x|_{I_\beta})^*(u_x|_{I_\beta})
=P(F_\beta,w).
\]
Let $V$ be the permutation matrix that puts the indices of each part
together. Use this order for both the rows and the columns. Then
\[
V^*PV=\bigoplus_\beta P[I_\beta,I_\beta]
     =\bigoplus_\beta P(F_\beta,w).
\]
Since $V^*PV=(V^*PV)^*=(V^*PV)^2$, each block $P(F_\beta,w)$ is an
orthogonal projection.

If $P=0$, formula \eqref{eq:gramsentries} gives $\sigma=0$. Thus $w=0$
and the claim holds. Suppose $P\ne0$, so $\sigma>0$. The components $C$
of the row graph on $F_{\beta_0}$ partition $\supp w$. Adding their
weights and using the hypothesis gives
\[
\sigma=\sum_Cw(C)\in p^{-d}\Zp.
\]
If the row graph on $F$ is connected, its only component has weight
$\sigma$, and the claim follows. This includes $n=1$, for which
$\supp w=\{0\}$.

Suppose the row graph on $F$ has at least two components. Let $L$ be
the size of the smallest coset that contains $\supp w$. Then $L>1$,
and Lemma~\ref{lem:coset} gives $F[L]=0$. Use the coset reduction from
Section~2.2 and write $\widetilde P=P(\widetilde F,\widetilde w)$.
Set $D_{x_0,\beta}=D_{x_0}[I_\beta,I_\beta]$. Since $D_{x_0}$ is
diagonal, restricting \eqref{eq:fold} to $I_\beta$ gives
\[
u_{\varphi(t)}|_{I_\beta}=(v_t|_{I_\beta})D_{x_0,\beta}
\qquad(t\in\Z_L).
\]
Use this identity on the two parts in the Gram sum to get
\[
\begin{aligned}
P[I_\beta,I_\gamma]
 &=\sum_{t\in\Z_L}\widetilde w(t)
     (u_{\varphi(t)}|_{I_\beta})^*(u_{\varphi(t)}|_{I_\gamma})\\
 &=D_{x_0,\beta}^*\,\widetilde P[I_\beta,I_\gamma]\,D_{x_0,\gamma}.
\end{aligned}
\]
The diagonal factors are unitary. For $\beta\ne\gamma$, this gives
\[
\widetilde P[I_\beta,I_\gamma]
=D_{x_0,\beta}\,P[I_\beta,I_\gamma]\,D_{x_0,\gamma}^*=0.
\]
By Section~2.2, $\widetilde P$ is an orthogonal projection,
$\widetilde F[L]=F[L]=0$, and $\widetilde w(\Z_L)=\sigma$.
The same reduction on $F_{\beta_0}$ preserves the component weights
on that sublist, so $\Rp_d$ still holds there.
We now work in $\Z_L$ and use $F,w,u_x,P$ for the reduced objects.
Put $T=L/q$. The column indices stay the same, so we use the same
$V$ to group the parts of the reduced $P$.

By \eqref{eq:gramsentries} and the block condition,
\[
z_\xi=z_\eta
\quad\Longrightarrow\quad P_{\xi\eta}=\sigma>0
\quad\Longrightarrow\quad \xi,\eta\in I_\beta\ \text{for some }\beta .
\]
Choose a matching as in Lemma~\ref{lem:split}. Let $\mathcal G$
index the groups. Write $z_g$ for the base frequency of group $g$
and set $F'=(z_g)_{g\in\mathcal G}$. For any order on $\mathcal G$,
the lemma gives
\[
U=RHSD,\qquad U^*PU=\bigoplus_{s=0}^{q-1}Q_s,\qquad
Q_s=P(F',w_s),\qquad w_s(x)=q\,w(s+qx)\quad(x\in\Z_T).
\]
The $q$ columns of a matched group have distinct frequencies
and may belong to different parts of the partition.

\smallskip
\noindent\textit{Case 1. Every matched group lies in one part.}
Let $\mathcal G_\beta$ be the groups in $I_\beta$ and set
$F_\beta'=(z_g)_{g\in\mathcal G_\beta}$. The part $I_\beta$ is a union
of whole groups. Thus for each $0\le z<T$, the frequencies
$z,z+T,\ldots,z+(q-1)T$ occur equally often in $F_\beta$.
By Lemma~\ref{thm:coeff}, $F_\beta[L]=0$.

Order the groups in each $\mathcal G_\beta$. To order $\mathcal G$,
take these lists in the block order of $V^*PV$.
For each $\beta$, apply Lemma~\ref{lem:split} to $P(F_\beta,w)$
using the chosen order on $\mathcal G_\beta$. This gives
$W_\beta=R_\beta H_\beta S_\beta D_\beta$ and
\[
W_\beta^*P(F_\beta,w)W_\beta
=\bigoplus_{s=0}^{q-1}P(F_\beta',w_s).
\]
Applying $V$ and then $\bigoplus_\beta W_\beta$ puts the coordinates
in the order $(\beta,s,g)$. Let $\Pi$ change this to $(s,\beta,g)$.
The factors in $U$ then satisfy
\[
R=V\Bigl(\bigoplus_\beta R_\beta\Bigr),\qquad
H=\bigoplus_\beta H_\beta,\qquad
S=\Bigl(\bigoplus_\beta S_\beta\Bigr)\Pi .
\]
Both $D$ and the corresponding $D_\beta$ use the phase
$\zeta_L^{-sz_g}$ at $(s,g)$. Hence
$\Pi D=(\bigoplus_\beta D_\beta)\Pi$ and
$U=V(\bigoplus_\beta W_\beta)\Pi$. The block form given by $V$ now gives
\[
\begin{aligned}
U^*PU
 &=\Pi^*\Bigl(\bigoplus_\beta W_\beta^*P(F_\beta,w)W_\beta\Bigr)\Pi\\
 &=\bigoplus_{s=0}^{q-1}\bigoplus_\beta P(F_\beta',w_s).
\end{aligned}
\]
Each $Q_s$ is block diagonal with parts $\mathcal G_\beta$. The blocks satisfy
$
Q_s[\mathcal G_\beta,\mathcal G_\beta]=P(F_\beta',w_s).
$

By hypothesis, $w$ satisfies $\Rp_d$ on $F_{\beta_0}$. Apply
Lemma~\ref{lem:split} to $P(F_{\beta_0},w)$. It gives $\Rp_d$ for $w_s$
on $F_{\beta_0}'$ for every $s$. This sublist has
$|F_{\beta_0}'|=|I_{\beta_0}|/q>0$ entries. Since $T=L/q<n$, we can
apply the induction hypothesis to each $Q_s=P(F',w_s)$. Hence each
$w_s$ satisfies $\Rp_d$ on $F'$. Now apply Lemma~\ref{lem:split}
to the full list to get $\Rp_d$ for $w$ on $F$.

\smallskip
\noindent\textit{Case 2. Some matched group meets more than one part.}
Choose a group of this kind with base frequency $z_*$, and let
$I_{\beta_*}$ be one of the parts it meets. For each base frequency
$z$ in $F'$ and $0\le j<q$, put $\tau_{z,j}=1$ if the columns of
frequency $z+jT$ lie in $I_{\beta_*}$ and $0$ otherwise.
All copies of a given frequency lie in the same part, so this defines
$\tau_{z,j}$ unambiguously. For $z=z_*$, both values $0$ and $1$ occur
among the $\tau_{z_*,j}$.

Let $E=\diag_{\xi\in I}(\one_{I_{\beta_*}}(\xi))$ be the coordinate
projection onto $I_{\beta_*}$. Right multiplication by $E$ keeps the
columns in $I_{\beta_*}$. Left multiplication keeps the rows in that
part. Since $P$ is block diagonal, both products keep the block
$P[I_{\beta_*},I_{\beta_*}]$ and have zeros everywhere else. Thus $PE=EP$.
In group order,
\[
R^*ER=\bigoplus_{g\in\mathcal G}
                    \diag_{0\le j<q}(\tau_{z_g,j}).
\]
For a group $g$ with $z=z_g$, conjugate its diagonal block by $H_q$.
The $(s,t)$ entry is
\[
\theta_{st}(z)
=\sum_{j=0}^{q-1}\overline{(H_q)_{j,s}}\,\tau_{z,j}(H_q)_{j,t}
=\frac1q\sum_{j=0}^{q-1}\tau_{z,j}\zeta_q^{(s-t)j}.
\]
The permutation $S$ orders the indices by channel. Before conjugating
by $D$, we therefore have entry $\delta_{g,g'}\theta_{st}(z_g)$
at $(s,g),(t,g')$.
%The channel blocks of $B=U^*EU$ are thus
%\begin{equation}\label{eq:indicator}
%(B_{st})_{g,g'}
%=\zeta_L^{sz_g}\delta_{g,g'}\theta_{st}(z_g)\zeta_L^{-tz_{g'}}
%=\delta_{g,g'}\zeta_L^{(s-t)z_g}\theta_{st}(z_g).
%\end{equation}
Write $B=U^*EU=D^*(S^*H^*R^*ERHS)D$. To get $B_{st}$, take the
$(s,t)$ block of $S^*H^*R^*ERHS$. Multiply it on the left by the
$s$th block of $D^*$ and on the right by the $t$th block of $D$.
This gives a product of three diagonal matrices on $\mathcal G$.
\[
B_{st}
=\diag_{g\in\mathcal G}\bigl(\zeta_L^{s z_g}\bigr)\,
\diag_{g\in\mathcal G}\bigl(\theta_{st}(z_g)\bigr)\,
\diag_{g\in\mathcal G}\bigl(\zeta_L^{-t z_g}\bigr).
\]
Multiplying the corresponding diagonal entries gives
\begin{equation}\label{eq:indicator}
	(B_{st})_{g,g'}
	=\delta_{g,g'}\zeta_L^{(s-t)z_g}\theta_{st}(z_g).
\end{equation}

For $s\ne t$, the root $\zeta_q^{s-t}$ has order $q$.
The coefficients $\tau_{z_*,j}$ take two distinct values, so
Lemma~\ref{thm:coeff} gives $\theta_{st}(z_*)\ne0$.

Conjugate $PE=EP$ by $U$ to get $
\Bigl(\bigoplus_{r=0}^{q-1}Q_r\Bigr)B
=B\Bigl(\bigoplus_{r=0}^{q-1}Q_r\Bigr)$.
The $(s,t)$ block gives $Q_sB_{st}=B_{st}Q_t$. For distinct $s,t$
with $Q_s\ne0$, choose $x\in\supp w_s$. The character row
$u'_x=(\zeta_T^{xz_g})_{g\in\mathcal G}$ satisfies $u'_xQ_s=u'_x$
by Section~2.1. Hence $v=u'_xB_{st}$ satisfies
\[
vQ_t=u'_xB_{st}Q_t=u'_xQ_sB_{st}=u'_xB_{st}=v.
\]
Since $B_{st}$ is diagonal on the group indices, \eqref{eq:indicator} gives
\[
v_g=\zeta_T^{xz_g}\zeta_L^{(s-t)z_g}\theta_{st}(z_g)
    =\zeta_L^{(qx+s-t)z_g}\theta_{st}(z_g).
\]
An entry with base frequency $z_*$ is nonzero. Since $vQ_t=v$,
we must have $Q_t\ne0$. In Lemma~\ref{lem:digits} take $h=q$, $a_0=qx+s-t$
(using an integer representative of $x$), and $b_z=\theta_{st}(z)$.
Then $hT=L$, $q\nmid a_0$, and $b_z\in\Q(\zeta_q)$. Since $b_{z_*}\ne0$,
the lemma gives saturation of $z_*$ in $Q_t$.

Since $P\ne0$, some $Q_{s_0}$ is nonzero. Apply the argument above with
this starting channel. It gives saturation in every $Q_t$ with
$t\ne s_0$ and shows that these projections are nonzero. Choose one
such $t$ and repeat the argument with the two channels interchanged.
This gives saturation in $Q_{s_0}$ as well, so $z_*$ is saturated
in every channel.

Let $c_*$ be the multiplicity of $z_*$ in $F'$ and put
$\sigma_s=\sum_{x\in\Z_T}w_s(x)$. By \eqref{eq:saturation},
$c_*\sigma_s=1$ for every $s$. Summing and using the unique
representation $s+qx$ of points in $\Z_L$ gives
\[
q=\sum_{s=0}^{q-1}c_*\sigma_s
 =q c_*\sum_{s=0}^{q-1}\sum_{x\in\Z_T}w(s+qx)
 =q c_*\sigma.
\]
Each group with base frequency $z_*$ contributes one copy of that
frequency to the full reduced list $F$. Its multiplicity in $F$ is
therefore $c_*$. The identity $c_*\sigma=1$ shows that $z_*$ is
saturated in $P$. Since $\sigma\in p^{-d}\Zp$,
Lemma~\ref{lem:saturated} gives $\Rp_d$ on the full reduced list.

In either case, Section~2.2 transfers $\Rp_d$ to the original list.
\end{proof}

\subsection{Diagonal phases on the parts}

Partition $I$ into $I_0,\ldots,I_{p-1}$, allowing empty parts.
Write $F_\beta=(z_\xi)_{\xi\in I_\beta}$. Put the indices of each
part together and use the same order for the rows and columns of
all matrices. For a matrix $A$ on $I$, the block
$A[I_\beta,I_\gamma]$ has rows in $I_\beta$ and columns in $I_\gamma$.

Put $\omega=\zeta_p$. For $0\le i<p$, let $P_i=P(F,w_i)$ be an
orthogonal projection. Let $\Omega_i$ be the diagonal unitary
with entry $\omega^{-i\beta}$ at each $\xi\in I_\beta$. We assume
\begin{equation}\label{eq:twist}
\begin{aligned}
P_i\Omega_i^*\Omega_jP_j&=0 &&(i\ne j),\\
H_\beta&=\Bigl(\sum_{i=0}^{p-1}P_i\Bigr)[I_\beta,I_\beta]
&&\text{is an orthogonal projection for each }\beta .
\end{aligned}
\end{equation}
For a nonempty part $I_\beta$, restrict the character rows to that part.
Then
\[
H_\beta
=\sum_{x\in\Z_n}\Bigl(\sum_{i=0}^{p-1}w_i(x)\Bigr)
 (u_x|_{I_\beta})^*(u_x|_{I_\beta})
=P\Bigl(F_\beta,\sum_{i=0}^{p-1}w_i\Bigr).
\]

\begin{lemma}\label{lem:twist}
Assume \eqref{eq:twist}. For $\beta\ne\gamma$,
\[
P_0[I_\beta,I_\gamma]=\cdots=P_{p-1}[I_\beta,I_\gamma].
\]
If $x\in\supp w_i$ and $y\in\supp w_j$ with $i\ne j$, then
\begin{equation}\label{eq:maskequal}
M_{F_0}(\zeta_n^{x-y})=\cdots=M_{F_{p-1}}(\zeta_n^{x-y}).
\end{equation}
\end{lemma}

\begin{proof}
Let $S_i=\Omega_iP_i\Omega_i^*$ and $S=\sum_{i=0}^{p-1}S_i$.
Each $S_i$ is an orthogonal projection since $\Omega_i$ is unitary.
For $i\ne j$, \eqref{eq:twist} gives
\[
S_iS_j=\Omega_i(P_i\Omega_i^*\Omega_jP_j)\Omega_j^*=0.
\]
Thus $S^*=S$ and $S^2=\sum_i S_i^2+\sum_{i\ne j}S_iS_j=S$.

Left multiplication by $\Omega_i$ scales rows in $I_\beta$ by
$\omega^{-i\beta}$, and right multiplication by $\Omega_i^*$ scales
columns in $I_\gamma$ by $\omega^{i\gamma}$. Hence
\[
S[I_\beta,I_\gamma]
=\sum_{i=0}^{p-1}\omega^{-i(\beta-\gamma)}P_i[I_\beta,I_\gamma].
\]
On the diagonal block, the phases cancel and
$S[I_\beta,I_\beta]=H_\beta$. The identity $S^2=S$ and the relation
$S[I_\gamma,I_\beta]=S[I_\beta,I_\gamma]^*$ give
\[
\begin{aligned}
H_\beta
&=(S^2)[I_\beta,I_\beta]
 =\sum_{\gamma=0}^{p-1}S[I_\beta,I_\gamma]S[I_\gamma,I_\beta]\\
&=H_\beta^2+\sum_{\gamma\ne\beta}
 S[I_\beta,I_\gamma]S[I_\beta,I_\gamma]^*.
\end{aligned}
\]
Since $H_\beta^2=H_\beta$, take the diagonal entry at $\xi\in I_\beta$
to get
\[
0=\sum_{\gamma\ne\beta}\sum_{\eta\in I_\gamma}|S_{\xi\eta}|^2.
\]
All the terms are nonnegative. Thus $S_{\xi\eta}=0$ whenever $\xi$
and $\eta$ are in different parts.

Put $K=\Q(\zeta_n)$. Since $p$ and $n$ are coprime,
$K(\omega)=\Q(\zeta_{pn})$. Thus
\[
[K(\omega):K]=\frac{\varphi(pn)}{\varphi(n)}=p-1.
\]
Since $\Phi_p(\omega)=0$ and $\deg\Phi_p=p-1$,
$\Phi_p$ is the minimal polynomial of $\omega$ over $K$ and is
irreducible over $K$.

Fix $\xi\in I_\beta$ and $\eta\in I_\gamma$ with $\beta\ne\gamma$.
Since $S_{\xi\eta}$ is zero, we have
\[
\sum_{i=0}^{p-1}(P_i)_{\xi\eta}
 \bigl(\omega^{-(\beta-\gamma)}\bigr)^i=0.
\]
The coefficients belong to $K$ by \eqref{eq:gramsentries}.
Since $\beta\not\equiv\gamma\pmod p$, the root
$\omega^{-(\beta-\gamma)}$ has order $p$. Lemma~\ref{thm:coeff}
then says that $(P_i)_{\xi\eta}$ is the same for all $i$.
Hence the off-diagonal blocks agree.

For $x\in\supp w_i$ and $y\in\supp w_j$ with $i\ne j$,
Section~2.1 gives $u_xP_i=u_x$ and $u_yP_j=u_y$. Take the
adjoint of $u_yP_j=u_y$ and use $P_j^*=P_j$ to get
$P_ju_y^*=u_y^*$. Then
\[
0=u_x(P_i\Omega_i^*\Omega_jP_j)u_y^*
 =(u_xP_i)\Omega_i^*\Omega_j(P_ju_y^*)
 =u_x\Omega_i^*\Omega_j u_y^*.
\]
On $I_\beta$, the diagonal entries of $\Omega_i^*\Omega_j$ are
$\omega^{(i-j)\beta}$. Expanding the inner product gives
\[
0
=\sum_{\beta=0}^{p-1}\omega^{(i-j)\beta}
  \sum_{\xi\in I_\beta}\zeta_n^{(x-y)z_\xi}
=\sum_{\beta=0}^{p-1}\omega^{(i-j)\beta}
  M_{F_\beta}(\zeta_n^{x-y}).
\]
Each mask value lies in $K$. The root $\omega^{i-j}$ has order $p$,
so Lemma~\ref{thm:coeff} gives \eqref{eq:maskequal}.
\end{proof}

\subsection{The loss of one factor of \texorpdfstring{$p$}{p}}

\Needspace{5\baselineskip}
\begin{theorem}\label{thm:weights}
Assume \eqref{eq:twist}, and let $d\ge0$ be an integer. Suppose
$H_{\beta_0}\ne0$ for some $\beta_0$ and $\sum_i w_i$ satisfies
$\Rp_d$ on $F_{\beta_0}$. Then each $w_i$ satisfies $\Rp_{d+1}$ on $F$.
\end{theorem}

\begin{proof}
Fix $d$. We use strong induction on the exponent of $q$ in $n$.
Put $w=\sum_i w_i$ and $\sigma_i=w_i(\Z_n)$. Take the row graph on
$F_{\beta_0}$ with weights $w$. By adding the weights of its
components $C$, we get
\begin{equation}\label{eq:summass}
\sum_{i=0}^{p-1}\sigma_i=\sum_C w(C)\in p^{-d}\Zp.
\end{equation}

\smallskip
\noindent\textit{A zero projection.}
Suppose one of the $P_i$ is zero. By Lemma~\ref{lem:twist}, all the
projections have the same off-diagonal blocks. Those blocks are
zero, so every $P_i$ is block diagonal. Each compression
$P(F_{\beta_0},w_i)$ is therefore an orthogonal projection.
On $I_{\beta_0}$, the phase $\Omega_i^*\Omega_j$ acts by the scalar
$\omega^{(i-j)\beta_0}$. Take the corresponding diagonal block
of the first relation in \eqref{eq:twist} to get
\[
0=\omega^{(i-j)\beta_0}P(F_{\beta_0},w_i)P(F_{\beta_0},w_j)
\qquad(i\ne j).
\]
These compressions are pairwise orthogonal and sum to $H_{\beta_0}$.
Lemma~\ref{lem:orthogonal-sums} gives $\Rp_d$ for each $w_i$ on
$F_{\beta_0}$. Apply Lemma~\ref{lem:parts} to get $\Rp_d$ on $F$
for each $w_i$. It follows that $\Rp_{d+1}$ holds as well, since
$p^{-d}\Zp\subseteq p^{-(d+1)}\Zp$.

\smallskip
\noindent\textit{The common support coset.}
Assume all $P_i$ are nonzero. Then each $\supp w_i$ is nonempty.
Fix an index $i_0$ and choose $x_0\in\supp w_{i_0}$.
Let $R$ be the size of the smallest coset that contains
$Y=\bigcup_i\supp w_i=\supp w$. This coset is $x_0+(n/R)\Z_n$.
Apply the reduction of Section~2.2 to all the weights using
$\varphi(t)=x_0+(n/R)t$ for $t\in\Z_R$.
Write $\widetilde F$ and $\widetilde F_\beta$ for the reduced lists.
Set $\widetilde w_i=w_i\circ\varphi$ and
$\widetilde w=\sum_i\widetilde w_i=w\circ\varphi$.
By \eqref{eq:fold}, $\widetilde P_i=D_{x_0}P_iD_{x_0}^*$.
The diagonal unitary $D_{x_0}$ is the same for all the projections.
Since $D_{x_0}$ commutes with each $\Omega_i$, we have
\[
\widetilde P_i\Omega_i^*\Omega_j\widetilde P_j
=D_{x_0}(P_i\Omega_i^*\Omega_jP_j)D_{x_0}^*=0
\qquad(i\ne j).
\]
For a nonempty part $I_\beta$, put
$D_{x_0,\beta}=D_{x_0}[I_\beta,I_\beta]$. Then
\[
\widetilde H_\beta
:=\Bigl(\sum_i\widetilde P_i\Bigr)[I_\beta,I_\beta]
=D_{x_0,\beta}H_\beta D_{x_0,\beta}^*
=P(\widetilde F_\beta,\widetilde w).
\]
Since $D_{x_0,\beta}$ is unitary, $\widetilde H_\beta$ is an
orthogonal projection. We also have $\widetilde H_{\beta_0}\ne0$.
By Section~2.2, $\widetilde w$ satisfies $\Rp_d$ on
$\widetilde F_{\beta_0}$. Each $\widetilde P_i$ has the same component
weights as $P_i$. In particular,
$\widetilde w_i(\Z_R)=\sigma_i$.

For $r\mid R$ and $c\in\Z_R$, we have
\[
\varphi\left(c+\frac{R}{r}k\right)
=x_0+\frac{n}{R}c+\frac{n}{r}k
\quad\text{in }\Z_n,\qquad 0\le k<r.
\]
Thus a coset of size $r$ in $\Z_R$ maps onto a coset of size $r$
in $\Z_n$. If a coset of size $r<R$ contained the union of the
reduced supports, its image would contain $Y$. This would contradict
the choice of $R$. Therefore the smallest coset that contains the union
of the reduced supports has size $R$.
Since $\varphi(0)=x_0$, we have
$\widetilde w_{i_0}(0)=w_{i_0}(x_0)>0$. We now work in $\Z_R$
and drop the tildes, so $0\in\supp w_{i_0}$.

If $R=1$, then $\supp w_i=\{0\}$ and $\sigma_i=w_i(0)$ for every $i$.
Put $u_0=(1,\ldots,1)$ and $f=|I|$. We have
$P_i=\sigma_i u_0^*u_0$ and $u_0u_0^*=f$. Squaring gives
$
P_i^2=\sigma_i^2 u_0^*(u_0u_0^*)u_0=f\sigma_i P_i.
$
Since $P_i^2=P_i\ne0$, we get $\sigma_i=1/f$ for every $i$. By
\eqref{eq:summass},
\[
\sigma_i=\frac{1}{p}\sum_{j=0}^{p-1}\sigma_j
\in p^{-(d+1)}\Zp.
\]
The only component of each row graph is $\{0\}$, with weight
$\sigma_i$. Hence each $w_i$ satisfies $\Rp_{d+1}$ on $F$.
Together with the zero-projection case, this covers $n=1$.

Suppose $R>1$. Lemma~\ref{lem:coset} gives $x\in\supp w_i$ and
$y\in\supp w_j$ for some $i\ne j$ with $\ord(x-y)=R$.
Then $\zeta_R^{x-y}$ is a primitive $R$th root. By
Lemma~\ref{lem:twist}, all the part masks take the same value at
this root. Their coefficients are integers, so conjugacy of the
primitive roots gives
\begin{equation}\label{eq:diffmasks}
(M_{F_\beta}-M_{F_\gamma})(\zeta_R)=0
\qquad(0\le\beta,\gamma<p).
\end{equation}

Put $T=R/q$. For $0\le z<T$ and $0\le j<q$, let
\[
n_{j,\beta}(z)=|\{\xi\in I_\beta:z_\xi=z+jT\}|,\qquad
\boldsymbol n_j(z)=(n_{j,0}(z),\ldots,n_{j,p-1}(z)).
\]
Apply Lemma~\ref{thm:coeff} to the differences in \eqref{eq:diffmasks}.
It follows that $n_{j,\beta}(z)-n_{j,\gamma}(z)$ is independent of $j$.
Hence $n_{j,\beta}(z)-n_{j',\beta}(z)=n_{j,\gamma}(z)-n_{j',\gamma}(z)$
for all $\beta,\gamma$. Equivalently,
\begin{equation}\label{eq:countdiff}
\boldsymbol n_j(z)-\boldsymbol n_{j'}(z)\in\Z\,\one,
\end{equation}
where $\one$ is the all-one vector of length $p$.

\smallskip
\noindent\textit{Case 1. The common mask value at $\zeta_R$ is zero.}
Since $F_\beta[R]=0$ for each $\beta$, Lemma~\ref{thm:coeff} lets
us match the columns in each part into complete $q$-point groups.
Let $\mathcal G_\beta$ index the groups in $I_\beta$ and write
$z_g$ for the base frequency of group $g$. Set
$\mathcal G=\bigsqcup_\beta\mathcal G_\beta$,
$F'=(z_g)_{g\in\mathcal G}$, and $F_\beta'=(z_g)_{g\in\mathcal G_\beta}$.
Apply Lemma~\ref{lem:split} to every $P_i$ using this matching of $F$.
The lemma uses the same column change $U$ in each case and gives
\[
U^*P_iU=\bigoplus_{s=0}^{q-1}P_{i,s},\qquad
P_{i,s}=P(F',w_{i,s}),\qquad w_{i,s}(t)=q\,w_i(s+qt)
\quad(t\in\Z_T).
\]

Let $\Omega_i'$ be the diagonal matrix on $\mathcal G$ with
$(\Omega_i')_{gg}=\omega^{-i\beta}$ for $g\in\mathcal G_\beta$.
In group order, the block of $\Omega_i$ on group $g$ is
$(\Omega_i')_{gg}I_q$. With the column change from
Lemma~\ref{lem:split}, this gives
\[
\begin{aligned}
	U^*\Omega_iU
	&=D^*S^*\Bigl(\bigoplus_{g\in\mathcal G}
	(\Omega_i')_{gg}H_q^*H_q\Bigr)SD\\
	&=D^*S^*\Bigl(\bigoplus_{g\in\mathcal G}
	(\Omega_i')_{gg}I_q\Bigr)SD\\
	&=D^*\Bigl(\bigoplus_{s=0}^{q-1}\Omega_i'\Bigr)D
	=\bigoplus_{s=0}^{q-1}\Omega_i'.
\end{aligned}
\]
Here $S$ orders the group scalars by channel. The resulting matrix
is diagonal and commutes with $D$.
Conjugate the first relation in \eqref{eq:twist} by $U$ to get
\[
0=U^*(P_i\Omega_i^*\Omega_jP_j)U
=\bigoplus_{s=0}^{q-1}
\bigl(P_{i,s}(\Omega_i')^*\Omega_j'P_{j,s}\bigr)
\qquad(i\ne j).
\]
Thus $P_{i,s}(\Omega_i')^*\Omega_j'P_{j,s}=0$ in every channel $s$.

For each nonempty $I_\beta$, apply Lemma~\ref{lem:split} to
$H_\beta=P(F_\beta,w)$ with the groups in $\mathcal G_\beta$.
The channel blocks are the orthogonal projections
\[
H_{\beta,s}=P(F_\beta',w_s),\qquad
w_s(t)=q\,w(s+qt)=\sum_{i=0}^{p-1}w_{i,s}(t).
\]

Let $u'_t=(\zeta_T^{t z_g})_{g\in\mathcal G}$ be the character row
on $F'$. Restrict this row to $\mathcal G_\beta$ to get the
character row on $F_\beta'$. Substituting the formula for $w_s$ gives
\[
\begin{aligned}
	H_{\beta,s}
	&=\sum_{t\in\Z_T}w_s(t)
	(u'_t|_{\mathcal G_\beta})^*(u'_t|_{\mathcal G_\beta})
	=\sum_{i=0}^{p-1}\sum_{t\in\Z_T}w_{i,s}(t)
	(u'_t|_{\mathcal G_\beta})^*(u'_t|_{\mathcal G_\beta})\\
	&=\sum_{i=0}^{p-1}P_{i,s}[\mathcal G_\beta,\mathcal G_\beta]
	=\Bigl(\sum_{i=0}^{p-1}P_{i,s}\Bigr)[\mathcal G_\beta,\mathcal G_\beta].
\end{aligned}
\]
For an empty part, the block is the zero-dimensional projection.
Thus the second condition of \eqref{eq:twist} holds in each channel.

The sublist $F_{\beta_0}'$ has $|I_{\beta_0}|/q>0$ entries.
Apply Lemma~\ref{lem:split} to $H_{\beta_0}$. The correspondence
of components gives $\Rp_d$ for $w_s$ on $F_{\beta_0}'$ for
every $s$. Fix a channel $s$. If all $P_{i,s}$ are zero, the
diagonal formula gives zero weights $w_{i,s}$, so $\Rp_{d+1}$ holds.
Otherwise some $w_{i,s}$ has positive total weight. Each diagonal
entry of $H_{\beta_0,s}$ is then $\sum_i w_{i,s}(\Z_T)>0$, so
$H_{\beta_0,s}\ne0$. Since $T=R/q<n$, apply the induction hypothesis
at the same $d$ to the family $(P_{i,s})_i$ on $\Z_T$. It gives
$\Rp_{d+1}$ for every $w_{i,s}$ on $F'$. Lemma~\ref{lem:split}
now gives $\Rp_{d+1}$ for every $w_i$ on the reduced list $F$.
Use Section~2.2 to return to the original list.

\smallskip
\noindent\textit{Case 2. The common mask value at $\zeta_R$ is nonzero.}
By Lemma~\ref{thm:coeff}, there is a cycle with two positions for which
$\boldsymbol n_j(z)\ne\boldsymbol n_{j'}(z)$. By \eqref{eq:countdiff},
we can order the positions so that
$\boldsymbol n_j(z)=\boldsymbol n_{j'}(z)+c\one$ with $c>0$.
Thus $n_{j,\beta}(z)\ge c$ for all $\beta$ and the frequency
$z+jT$ occurs in every part. Choose two columns $\xi,\eta$ at this
frequency in different parts. Formula \eqref{eq:gramsentries} gives
$(P_i)_{\xi\eta}=\sigma_i$. Lemma~\ref{lem:twist} says that these
entries agree for all $i$. Write $\sigma$ for their common value.
Then \eqref{eq:summass} gives
\begin{equation}\label{eq:massbound}
p\sigma=\sum_i\sigma_i\in p^{-d}\Zp,\qquad
\sigma\in p^{-(d+1)}\Zp .
\end{equation}
If the row graph of $P_{i_0}$ is connected, its only component
has weight $\sigma$. Thus \eqref{eq:massbound} gives $\Rp_{d+1}$,
and Section~2.2 gives the same conclusion in the original group.
Suppose now that this graph is disconnected.

Let $L$ be the size of the smallest coset containing
$\supp w_{i_0}$. By Lemma~\ref{lem:coset}, $L>1$ and $F[L]=0$.
By \eqref{eq:diffmasks}, $F[R]=pF_\beta[R]\ne0$. Hence $L<R$.
Both $L$ and $R$ are powers of $q$, so $L\mid T=R/q$.
Since $0\in\supp w_{i_0}$, we have
\[
\supp w_{i_0}\subseteq (R/L)\Z_R\subseteq q\Z_R.
\]

For $0\le z<L$, let $C_z=\{\xi\in I:z_\xi\equiv z\pmod L\}$.
Call $C_z$ mixed if it has indices in at least two parts. Let $J$
be the union of the mixed classes and put $E=I\setminus J$. The
indices of the frequency found above belong to a mixed class,
so $J\ne\varnothing$.

Take $\xi,\eta$ in the same mixed class and in different parts.
Their frequency difference is a multiple of $L$. Thus
$x(z_\eta-z_\xi)\equiv0\pmod R$ for $x\in\supp w_{i_0}$.
Formula \eqref{eq:gramsentries} gives $(P_{i_0})_{\xi\eta}=\sigma$.
By Lemma~\ref{lem:twist}, $(P_i)_{\xi\eta}=\sigma$ for all $i$.
Both diagonal entries equal $\sigma$. We also have
$(P_i)_{\eta\xi}=\overline{(P_i)_{\xi\eta}}=\sigma$. Hence
\[
\begin{aligned}
\|P_i(\mathbf e_\xi-\mathbf e_\eta)\|^2
&=(\mathbf e_\xi-\mathbf e_\eta)^*P_i(\mathbf e_\xi-\mathbf e_\eta)\\
&=(P_i)_{\xi\xi}-(P_i)_{\xi\eta}-(P_i)_{\eta\xi}
   +(P_i)_{\eta\eta}=0,
\end{aligned}
\]
Here $\mathbf e_\xi,\mathbf e_\eta$ are the corresponding standard
columns. Thus $P_i\mathbf e_\xi=P_i\mathbf e_\eta$. To compare two
indices in the same part of a mixed class, compare both with one
index in another part of that class. This proves that the columns
of $P_i$ agree within each mixed class.

Fix $\xi\in J$ and a row index $\kappa\in I_\beta$. The mixed class
of $\xi$ meets at least two parts, so choose $\eta\notin I_\beta$
in that class. The columns at $\xi$ and $\eta$ agree in each $P_i$.
The entry at $(\kappa,\eta)$ is the same for all $i$ by
Lemma~\ref{lem:twist}. Hence
\[
(P_i)_{\kappa\xi}=(P_i)_{\kappa\eta}
=(P_{i_0})_{\kappa\eta}=(P_{i_0})_{\kappa\xi}.
\]
Since $\kappa\in I$ and $\xi\in J$ were arbitrary, we have
$P_i[I,J]=P_{i_0}[I,J]$ for every $i$.

\smallskip
\noindent\textit{A fixed row on $J$.}
The subgroup $q\Z_R$ has size $R/q<R$. By the choice of $R$, the
union of the reduced supports contains a point $y\notin q\Z_R$.
Then $y$ is a unit modulo $R$. Choose $j$ with $y\in\supp w_j$.
By Section~2.1, $u_yP_j=u_y$. If $\xi,\eta$ lie in the same mixed
class, the equality of their columns gives
\[
\zeta_R^{y z_\xi}-\zeta_R^{y z_\eta}
=u_yP_j(\mathbf e_\xi-\mathbf e_\eta)=0.
\]
As $y$ is a unit, all frequencies in a mixed class agree modulo $R$.

Write $F_E=(z_\xi)_{\xi\in E}$ and fix $0\le z<T$. Since $L\mid T$,
the frequencies $z+jT$, $0\le j<q$, all have the same residue
modulo $L$. Their column indices lie in a single class.
If the class meets at most one part, all its indices lie in $E$.
Choose a part $I_\gamma$ disjoint from the class. Then
$n_{j,\gamma}(z)=0$ for every $j$. By \eqref{eq:countdiff}, for
any $j,j'$ we have
\[
\boldsymbol n_j(z)-\boldsymbol n_{j'}(z)=c\one,\qquad
c=n_{j,\gamma}(z)-n_{j',\gamma}(z)=0.
\]
Thus the count vectors are equal. The multiplicity of $z+jT$ in
$F_E$ is $\sum_\beta n_{j,\beta}(z)$ and is independent of $j$.
If the class is mixed, all its indices lie in $J$. In that case
every $z+jT$ has multiplicity zero in $F_E$. In both cases the
multiplicities in $F_E$ are constant as $j$ varies. By
Lemma~\ref{thm:coeff}, $M_{F_E}(\zeta_R)=0$.

For $x\in\supp w_{i_0}$, we have $q\mid x$ and $q\nmid y$.
Thus $\zeta_R^{y-x}$ is a primitive $R$th root. Conjugacy of the
primitive roots gives
\[
(u_y|_E)(u_x|_E)^*=M_{F_E}(\zeta_R^{y-x})=0.
\]
Now use the Gram formula for the block with rows in $E$ and
columns in $J$. We get
\[
(u_y|_E)P_{i_0}[E,J]
=\sum_{x\in\supp w_{i_0}}
 w_{i_0}(x)\bigl((u_y|_E)(u_x|_E)^*\bigr)(u_x|_J)=0 .
\]
Using $P_j[I,J]=P_{i_0}[I,J]$ and $u_yP_j=u_y$, we obtain
\[
\begin{aligned}
u_y|_J
&=(u_yP_j)|_J=(u_yP_{i_0})|_J\\
&=(u_y|_J)P_{i_0}[J,J]+(u_y|_E)P_{i_0}[E,J]\\
&=(u_y|_J)P_{i_0}[J,J].
\end{aligned}
\]
Let $v$ be the extension of $u_y|_J$ by zero on $E$. Since
$I_f-P_{i_0}$ is an orthogonal projection,
\[
\begin{aligned}
\|v(I_f-P_{i_0})\|^2
&=v(I_f-P_{i_0})v^*\\
&=(u_y|_J)(I_{|J|}-P_{i_0}[J,J])(u_y|_J)^*=0.
\end{aligned}
\]
Hence $vP_{i_0}=v$.

\smallskip
\noindent\textit{Saturation on the smaller group.}
Put $h=R/L>1$ and write $\bar z_\xi=z_\xi\bmod L$ with
$0\le\bar z_\xi<L$. Since $\supp w_{i_0}\subseteq h\Z_R$,
use $t\mapsto ht$ to reduce to $\Z_L$. For $t\in\Z_L$,
\[
u_{ht}
=(\zeta_R^{htz_\xi})_{\xi\in I}
=(\zeta_L^{t\bar z_\xi})_{\xi\in I}.
\]
The list $(\bar z_\xi)_{\xi\in I}$ with weights
$t\mapsto w_{i_0}(ht)$ therefore gives the same projection $P_{i_0}$.
This projection still fixes $v$.

Consider a mixed class of residue $0\le z<L$. Its columns have a
common frequency modulo $R$. Write this frequency as $z+Lk_z$,
with $0\le k_z<h$. For each $\xi$ in the class,
\[
v_\xi=(u_y)_\xi
=\zeta_R^{y(z+Lk_z)}
=\zeta_{hL}^{yz}\zeta_h^{yk_z}.
\]
The entries of $v$ on all other classes are zero. Apply
Lemma~\ref{lem:digits} with $T=L$ and $a_0=y$, using an integer
representative of $y$. Set $b_z=\zeta_h^{yk_z}$ on mixed classes
and $b_z=0$ elsewhere, so $b_z\in\Q(\zeta_h)$. Both $L$ and $h$
are powers of $q$, and $q\nmid y$. Since $J\ne\varnothing$, some
frequency on the list has $b_z\ne0$. The lemma gives a saturated
frequency for the projection on $\Z_L$.

The total weight is $\sigma$. By Lemma~\ref{lem:saturated}, each
component $C$ has weight $q^{-r_C}\sigma$ for some $r_C\in\Z_{\ge0}$.
By \eqref{eq:massbound} and $q^{-1}\in\Zp$, these weights belong to
$p^{-(d+1)}\Zp$. Section~2.2 preserves the component weights in both
reductions, so $w_{i_0}$ satisfies $\Rp_{d+1}$ on the original list.
Since $i_0$ was arbitrary, the theorem follows.
\end{proof}

\section{From the spectral matrix to character Gram projections}\label{sec:spectral}

\subsection{Coordinates and complete cycles}

Return to the spectral pair $(A,B)$ in $\Z_N$ and the integer $k$
from Section~\ref{sec:obstruction}. A complete $q$-cycle in $\Z_N$
has the form
\[
C_r=\left\{r+j\frac{N}{q}:0\le j<q\right\},
\qquad 0\le r<N/q.
\]
Let $S$ be the union of all such cycles contained in $B$, and put
$R_0=B\setminus S$. Distinct cycles are disjoint because they are
cosets of the same subgroup.

\begin{lemma}\label{lem:cycles}
The set $S$ is nonempty, and $R_0+N/p^k=R_0$ in $\Z_N$.
\end{lemma}

\begin{proof}
Put $N_h=N/p^h$ for $0\le h\le k, M_{C_r}(X)=X^r\sum_{j=0}^{q-1}X^{jN/q}.$
For $0\le h<k$, we have $\zeta_{N_h}^{N/q}=\zeta_q^{p^h}$, a
primitive $q$th root. Thus
$
M_{C_r}(\zeta_{N_h})
=\zeta_{N_h}^{r}\sum_{j=0}^{q-1}(\zeta_q^{p^h})^j=0.
$
The cycles in $S$ are disjoint, so $S[N_h]=0$. By
\eqref{eq:firstzeros}, $B[N_h]=0$, and hence
\[
R_0[N_h]=B[N_h]-S[N_h]=0\qquad(0\le h<k).
\]
By definition, $R_0$ contains no complete $q$-cycle.

\smallskip
Let $M$ have exactly the prime divisors $p$ and $q$, and suppose
$D\subseteq\Z_M$ satisfies $D[M]=0$ and contains no complete
$q$-cycle. For $\ell\in\{p,q\}$ put
$
\Psi_\ell(X)=\sum_{j=0}^{\ell-1}X^{jM/\ell}.
$
For an integer $t$, the geometric sum gives
\[
\Psi_\ell(\zeta_M^t)
=\sum_{j=0}^{\ell-1}\zeta_\ell^{jt}
=\begin{cases}
\ell, & \ell\mid t,\\
0, & \ell\nmid t .
\end{cases}
\]
At $\zeta_M^t$, the product $(p-\Psi_p)(q-\Psi_q)M_D$ vanishes
if $p\mid t$ or $q\mid t$ by the first two factors. In the
remaining case $\zeta_M^t$ is primitive, so $M_D(\zeta_M^t)=0$
by conjugacy of the primitive roots. Thus the product is
divisible by $X^M-1$. Expanding it gives
\[
M_D\equiv \frac1p\Psi_p M_D+\frac1q\Psi_q M_D
            -\frac1{pq}\Psi_p\Psi_q M_D
\pmod{X^M-1}.
\]

Let $f(x)=\one_D(x)$ for $x\in\Z_M$. Modulo $X^M-1$, the term $f(y)X^{y+iM/p}$
contributes to $X^x$ when $y+iM/p\equiv x\pmod M$.
Thus the coefficient of $X^x$ in the remainder of $\Psi_pM_D$ is
\[
\sum_{i=0}^{p-1} f(x-iM/p)=\sum_{i=0}^{p-1} f(x+iM/p).
\]
The other terms give their coefficients in the same way, so
\begin{equation}\label{eq:averages}
\begin{aligned}
f(x)
&=\frac1p\sum_{i=0}^{p-1} f(x+iM/p)
 +\frac1q\sum_{j=0}^{q-1} f(x+jM/q)\\
&\quad-\frac1{pq}\sum_{i=0}^{p-1}\sum_{j=0}^{q-1}
       f(x+iM/p+jM/q).
\end{aligned}
\end{equation}

Fix $x\in\Z_M$ and set $m_{ij}=f(x+iM/p+jM/q)$ for
$0\le i<p$ and $0\le j<q$. Apply \eqref{eq:averages} at the
point of entry $(i,j)$ and reindex the sums modulo $p$ and $q$
to obtain
\[
m_{ij}
=\frac1p\sum_{i'=0}^{p-1}m_{i'j}
 +\frac1q\sum_{j'=0}^{q-1}m_{ij'}
 -\frac1{pq}\sum_{i'=0}^{p-1}\sum_{j'=0}^{q-1}m_{i'j'}.
\]
The three terms average column $j$, row $i$, and the whole array,
respectively. The right side therefore has the form $a_i+b_j$.
Each row runs through a complete $q$-cycle, so the assumption
on $D$ gives
\[
0=\min_j m_{ij}=a_i+\min_j b_j
\]
for each $i$. Thus all $a_i$ are equal, and so $m_{00}=m_{10}$.
As $x$ was arbitrary, $D+M/p=D$.

\smallskip
Let $\pi:\Z_M\to\Z_{M/p}$ be reduction modulo $M/p$ and put
$D'=\pi(D)$. The period $M/p$ gives $D=\pi^{-1}(D')$.
For each $x'\in D'$ with $0\le x'<M/p$, the full fiber consists
of $x'+\ell M/p$ for $0\le\ell<p$. Summing over these fibers gives
\[
M_D(X)=M_{D'}(X)\sum_{\ell=0}^{p-1}X^{\ell M/p}.
\]
In particular, for $d\mid M/p$ we have $D[d]=pD'[d]$.

If $D'$ contained a complete $q$-cycle with base $x'$, its step
would be $M/(pq)$. The complete $q$-cycle in $\Z_M$ with base $x'$
and step $M/q$ projects onto this cycle, since
\[
\pi(x'+jM/q)=x'+pjM/(pq)\pmod{M/p}
\qquad(0\le j<q),
\]
and $pj\bmod q$ runs through all residues modulo $q$. As
$D=\pi^{-1}(D')$, that full cycle in $\Z_M$ would lie in $D$.
Hence $D'$ contains no complete $q$-cycle.

\smallskip
Starting with $R_0$, define $R_{h+1}$ to be the image of $R_h$
modulo $N_{h+1}$ for $0\le h<k$. We show by induction that $R_h$
contains no complete $q$-cycle and
\[
R_h[N_j]=0\qquad(h\le j<k).
\]
Both properties hold for $R_0$. 
At step $h<k$, suppose $R_h$ has no complete $q$-cycle and
$R_h[N_j]=0$ for $h\le j<k$. The zero
$R_h[N_h]=0$ gives the period $N_{h+1}$ by the one-step argument.
Thus $R_h$ is the full preimage of its reduction $R_{h+1}$ modulo
$N_{h+1}$. By the quotient argument, $R_{h+1}$ has no complete
$q$-cycle and
\[
R_{h+1}[N_j]=p^{-1}R_h[N_j]=0\qquad(h+1\le j<k).
\]
After $k$ steps, $x\in R_0$ if and only if $x\bmod N_k\in R_k$,
so $R_0+N/p^k=R_0$. If $S=\varnothing$, the period $N/p$ from
the first step makes $B=R_0$ a union of complete $p$-cycles,
contrary to property~(ii) of Section~\ref{sec:obstruction}.
\end{proof}

\Needspace{6\baselineskip}

Choose integers $r_p,r_q$ with $r_pq^b\equiv1\pmod{p^a}$
and $r_qp^a\equiv1\pmod{q^b}$. Use ordinary Chinese remainder
coordinates for a row point $x$ and dual coordinates for a
column frequency $y$,
\[
x_p=x\bmod p^a,\qquad x_q=x\bmod q^b,\qquad
y_p=r_py\bmod p^a,\qquad y_q=r_qy\bmod q^b.
\]
The inverse formulas are
\[
x\equiv x_p r_pq^b+x_q r_qp^a\pmod N,\qquad
y\equiv q^b y_p+p^a y_q\pmod N.
\]
Substitution of the formula for $x$ into the character gives
\[
\chi_y(x):=\zeta_N^{xy}
=\zeta_{p^a}^{x_p r_p y}\zeta_{q^b}^{x_q r_q y}
=\zeta_{p^a}^{x_py_p}\zeta_{q^b}^{x_qy_q}.
\]

Let $c=|A|=|B|$ and $U_{x,y}=c^{-1/2}\chi_y(x)$ for $x\in A$
and $y\in B$. Both coordinate maps are additive isomorphisms.
In the dual coordinates the two relevant shifts are
\[
\begin{aligned}
N/q&\longmapsto
  \bigl(r_p p^a q^{b-1}\bmod p^a,\,
        r_q p^a q^{b-1}\bmod q^b\bigr)=(0,q^{b-1}),\\
N/p^k&\longmapsto
  \bigl(r_p p^{a-k}q^b\bmod p^a,\,
        r_q p^{a-k}q^b\bmod q^b\bigr)=(p^{a-k},0).
\end{aligned}
\]
A complete $q$-cycle through $y$ in $\Z_N$ is
$
\{y+jN/q:0\le j<q\}.
$
Write $(y_p,y_q)$ for the dual coordinates of $y$. Since $N/q$
has dual coordinates $(0,q^{b-1})$, the cycle is represented by
$
\{(y_p,y_q+jq^{b-1}):0\le j<q\}.
$
Let $z$ be the least nonnegative residue of $y_q$ modulo $q^{b-1}$, we can start the same cycle at
$(y_p,z)$ and write it as
\[
\{(y_p,z+jq^{b-1}):0\le j<q\},
\qquad 0\le z<q^{b-1}.
\]
Lemma~\ref{lem:cycles} gives $R_0+N/p^k=R_0$ in $\Z_N$.
The dual coordinate map sends the shift $N/p^k$ to
$(p^{a-k},0)$, so this equality becomes
$R_0+(p^{a-k},0)=R_0$.

\subsection{Row classes and commuting column projections}

For $0\le v\le k$, define the row and column classes
\[
\begin{aligned}
C_{v,\epsilon}
  &=\{x\in A:x_p\equiv\epsilon\pmod{p^v}\},
  &0\le\epsilon<p^v,\\
A_\alpha
  &=\{x\in A:x_q\equiv\alpha\pmod q\},
  &0\le\alpha<q,\\
D_{v,u}
  &=\{y\in B:y_p\equiv u\pmod{p^{a-v}}\},
  &0\le u<p^{a-v}.
\end{aligned}
\]
Put $X_{v,\epsilon,\alpha}=C_{v,\epsilon}\cap A_\alpha$. Each
$x\in A$ determines one pair $(\epsilon,\alpha)$ at level $v$,
so these intersections are disjoint and cover $A$.

Fix $1\le j\le k$. For $x\in C_{j,s}$, we have
$x\equiv x_p\equiv s\pmod{p^j}$. Grouping the terms by $s$ and
writing each $s=r+i p^{j-1}$ with $0\le r<p^{j-1}$ and $0\le i<p$,
we obtain
\[
\begin{aligned}
	0=A[p^j]
	&=M_A(\zeta_{p^j})
	=\sum_{x\in A}\zeta_{p^j}^{x}\\
	&=\sum_{s=0}^{p^j-1}\sum_{x\in C_{j,s}}\zeta_{p^j}^{s}
	=\sum_{s=0}^{p^j-1}|C_{j,s}|\,\zeta_{p^j}^{s}\\
	&=\sum_{r=0}^{p^{j-1}-1}\sum_{i=0}^{p-1}
	|C_{j,r+i p^{j-1}}|\,\zeta_{p^j}^{r+i p^{j-1}}.
\end{aligned}
\]
Lemma~\ref{thm:coeff} applied to the last sum shows that, for
each fixed $r$, the $p$ counts $|C_{j,r+i p^{j-1}}|$ are equal.
Since $C_{j-1,r}=\bigsqcup_{i=0}^{p-1}C_{j,r+i p^{j-1}}$,
each child has size
\[
|C_{j,r+i p^{j-1}}|
=\frac{|C_{j-1,r}|}{p}
\qquad(0\le i<p).
\]
 Start with $C_{0,0}=A$ and iterate.
The same coefficient comparison applied to
$A[q]=\sum_{\alpha=0}^{q-1}|A_\alpha|\zeta_q^\alpha=0$ gives
equal counts $|A_\alpha|$ that sum to $c$. We obtain
\begin{equation}\label{eq:sizes}
|C_{v,\epsilon}|=\frac{c}{p^v},\qquad
|A_\alpha|=\frac{c}{q}=:\mu .
\end{equation}

For $T\subseteq B$, let $U_T$ consist of the columns of $U$
indexed by $T$. These columns are orthonormal, so $U_T^*U_T=I_T$
and $\Pi_T=U_TU_T^*$ is the orthogonal projection in $\C^A$
onto their span. Its entries are
\begin{equation}\label{eq:kernel}
(\Pi_T)_{x,x'}
=\sum_{y\in T}U_{x,y}\overline{U_{x',y}}
=\frac1c\sum_{y\in T}\chi_y(x-x').
\end{equation}
For $X\subseteq A$, let $E_X$ be the diagonal projection on
$\C^A$ with diagonal entries $\one_X(x)$. Since
\[
(E_X\Pi_T-\Pi_T E_X)_{x,x'}
=\bigl(\one_X(x)-\one_X(x')\bigr)(\Pi_T)_{x,x'},
\]
$E_X$ commutes with $\Pi_T$ exactly when the entries between
$X$ and $A\setminus X$ vanish in both directions.

\Needspace{7\baselineskip}
\begin{lemma}\label{prop:row-class-projections}
Fix $v,u$, set $D=S\cap D_{v,u}$, and put $K_Y=U_D^*E_YU_D$ on
$\C^D$ for $Y\subseteq A$. For $X=X_{v,\epsilon,\alpha}$ and
$X'=X_{v,\epsilon',\alpha'}$, we have $\Pi_D E_X=E_X\Pi_D$ and
\begin{equation}\label{eq:gramproduct}
K_XK_{X'}=K_{X\cap X'}.
\end{equation}
In particular, the $K_{X_{v,\epsilon,\alpha}}$ are pairwise
orthogonal projections with sum $I_D$.
\end{lemma}

\begin{proof}
Take $x,x'\in A$ with $x_p\not\equiv x_p'\pmod{p^v}$. Put
$\delta=x-x'$ and choose an integer representative
$\delta_p=p^e d_0$, where $0\le e<v$ and $p\nmid d_0$.
Taking $d_0^{-1}$ modulo $p^a$, set
\[
\gamma_p\equiv 1+p^{v-e}t\,d_0^{-1}\pmod{p^a}
\qquad(0\le t<p^{a-v}).
\]
Then $\gamma_p\equiv1\pmod p$ and
$(\gamma_p-1)\delta_p\equiv tp^v\pmod{p^a}$.
Let $d$ be the order of $\zeta_N^\delta$. Row orthogonality gives
$M_B(\zeta_N^\delta)=0$. By the Chinese remainder theorem, choose an integer $\gamma$ such that
$\gamma\equiv\gamma_p\pmod{p^a}$ and $\gamma\equiv1\pmod{q^b}$. Since $\gcd(\gamma,N)=1$,
$\zeta_N^{\gamma\delta}$ is also a primitive $d$th root. Hence
\[
0=M_B(\zeta_N^{\gamma\delta})
=\sum_{y\in B}\chi_y(\gamma\delta).
\]
The coordinates of $\gamma\delta$ are $(\gamma_p\delta_p,\delta_q)$.
Using $(\gamma_p-1)\delta_p\equiv tp^v\pmod{p^a}$, we have
\[
\begin{aligned}
	\chi_y(\gamma\delta)
	&=\zeta_{p^a}^{\gamma_p\delta_p y_p}\zeta_{q^b}^{\delta_q y_q}\\
	&=\chi_y(\delta)\zeta_{p^a}^{tp^v y_p}
	=\chi_y(\delta)\zeta_{p^{a-v}}^{t y_p}.
\end{aligned}
\]
Substitute this expression into the zero sum, multiply by
$\zeta_{p^{a-v}}^{-tu}$, and average over $0\le t<p^{a-v}$. This gives
\[
0=\sum_{y\in B}\chi_y(\delta)
\left(\frac{1}{p^{a-v}}
\sum_{t=0}^{p^{a-v}-1}\zeta_{p^{a-v}}^{t(y_p-u)}\right)
=\sum_{y\in D_{v,u}}\chi_y(\delta).
\]
The average in parentheses is $1$ for $y\in D_{v,u}$ and $0$
otherwise. By \eqref{eq:kernel}, the entries of $\Pi_{D_{v,u}}$
between distinct classes modulo $p^v$ vanish, so $\Pi_{D_{v,u}}$
commutes with every $E_{C_{v,\epsilon}}$. For $v=0$, we have
$E_{C_{0,0}}=I_A$.

Whenever $T+\tau=T$ in the dual coordinates, summing over
$y+\tau$ in \eqref{eq:kernel} gives
\[
(\Pi_T)_{x,x'}
=\frac1c\sum_{y\in T}\chi_{y+\tau}(x-x')
=\chi_\tau(x-x')(\Pi_T)_{x,x'} .
\]
For $T=R_0$, Lemma~\ref{lem:cycles} gives the shift
$\tau=(p^{a-k},0)$ and the factor
$\chi_\tau(x-x')=\zeta_{p^k}^{x_p-x_p'}$. When
$x_p\not\equiv x_p'\pmod{p^v}$ and $v\le k$, this factor differs
from $1$, which forces $(\Pi_{R_0})_{x,x'}=0$. Thus $\Pi_{R_0}$
commutes with $E_{C_{v,\epsilon}}$.
Since $B=S\sqcup R_0$ and $U$ is unitary, we have
\[
\Pi_S+\Pi_{R_0}
=U_SU_S^*+U_{R_0}U_{R_0}^*
=UU^*
=I_A.
\]
Since $\Pi_{R_0}$ commutes with $E_{C_{v,\epsilon}}$, so does
$\Pi_S=I_A-\Pi_{R_0}$.

Write $h_y$ for the $y$th column of $U$. Since
$h_y^*h_z=\delta_{y,z}$ and $\Pi_T=\sum_{y\in T}h_yh_y^*$,
\[
\Pi_S\Pi_{D_{v,u}}
=\sum_{y\in S}\sum_{z\in D_{v,u}}h_y(h_y^*h_z)h_z^*
=\sum_{y\in D}h_yh_y^*
=\Pi_D.
\]
Both factors commute with $E_{C_{v,\epsilon}}$, so their product
does as well. For any $(y_p,y_q)\in D=S\cap D_{v,u}$, the points
$(y_p,y_q+jq^{b-1})$, $0\le j<q$, form a complete $q$-cycle
in $S$. Each point has the same first coordinate $y_p$, which
satisfies $y_p\equiv u\pmod{p^{a-v}}$. Hence the entire cycle
lies in $D$. Adding $(0,q^{b-1})$ takes each such cycle onto
itself, so $D+(0,q^{b-1})=D$, and the shift formula gives
\[
(\Pi_D)_{x,x'}=\zeta_q^{x_q-x_q'}(\Pi_D)_{x,x'}.
\]
It follows that $(\Pi_D)_{x,x'}=0$ when
$x_q\not\equiv x_q'\pmod q$, so $\Pi_D$ commutes with
$E_{A_\alpha}$. Since it commutes with both factors in
$E_X=E_{C_{v,\epsilon}}E_{A_\alpha}$, it commutes with $E_X$.

Using $\Pi_D E_{X'}=E_{X'}\Pi_D$ and
$\Pi_D U_D=U_D(U_D^*U_D)=U_D$, we obtain
\[
\begin{aligned}
K_XK_{X'}
&=U_D^*E_X\Pi_D E_{X'}U_D\\
&=U_D^*E_XE_{X'}\Pi_D U_D\\
&=U_D^*E_{X\cap X'}U_D=K_{X\cap X'}.
\end{aligned}
\]
We have $K_X^*=K_X$ since $E_X^*=E_X$. Taking $X'=X$
gives $K_X^2=K_X$, and $X\cap X'=\varnothing$ gives
$K_XK_{X'}=0$. The row classes cover $A$ disjointly, so
\[
\sum_{\epsilon,\alpha}K_{X_{v,\epsilon,\alpha}}
=U_D^*\left(\sum_{\epsilon,\alpha}E_{X_{v,\epsilon,\alpha}}\right)U_D
=U_D^*U_D=I_D .\qedhere
\]
\end{proof}

\subsection{One representative from each cycle}

Put $n=q^{b-1}$. Fix an order of the complete $q$-cycles in $S$,
and let $\xi$ index these cycles. In cycle $\xi$, choose the unique
representative $y_\xi=(y_{p,\xi},z_\xi)$ with $0\le z_\xi<n$.
The column at position $j$ in that cycle is
$(y_{p,\xi},z_\xi+jn)$ for $0\le j<q$.

Fix $v,u$ and set $D=S\cap D_{v,u}$. Let $I_{v,u}$ be the
indices of the cycles contained in $D$, in the fixed order. Put
$F_{v,u}=(z_\xi)_{\xi\in I_{v,u}}$ in $\Z_n$. Two cycles with the
same value of $z_\xi$ give two distinct entries of this list.

For $x\in A_\alpha$, let $t(x)\in\Z_n$ be the unique residue satisfying
$x_q=\alpha+q\,t(x)$ in $\Z_{q^b}$. It exists because
$x_q\equiv\alpha\pmod q$, and division by $q$ determines $t(x)$
modulo $n=q^b/q$. For $X=X_{v,\epsilon,\alpha}$, put
\begin{equation}\label{eq:spectralweights}
w_X(t)=\frac{1}{\mu}\bigl|\{x\in X:t(x)=t\}\bigr|,
\qquad t\in\Z_n ,
\end{equation}
where $\mu=c/q$ is from \eqref{eq:sizes}. Each $x\in X$ is
counted once in the sum over $t$, so
$w_X(\Z_n)=|X|/\mu$.

\Needspace{5\baselineskip}
\begin{lemma}\label{lem:grams}
If $F_{v,u}$ is nonempty, then $P(F_{v,u},w_X)$ is an orthogonal
projection for every $X=X_{v,\epsilon,\alpha}$. When $v=0$, we have
$P(F_{0,u},w_{A_\alpha})=I_{F_{0,u}}$ for every $\alpha$,
where $I_{F_{0,u}}$ is the identity on the list coordinates.
\end{lemma}

\begin{proof}
For each $\xi\in I_{v,u}$, order its columns as
$(\xi,0),\ldots,(\xi,q-1)$ and right-multiply them by the Fourier
matrix $q^{-1/2}(\zeta_q^{-js})_{0\le j,s<q}$. Label the $s$th
output from cycle $\xi$ by $(s,\xi)$.
For each $s=0,\ldots,q-1$, put the columns $(s,\xi)$ together
in the fixed order of $\xi$. Write $U_DV$ for the resulting matrix.
The Fourier matrices and the column permutations are unitary,
so $V$ is unitary.

For $x\in A_\alpha$, the column at $(\xi,j)$ has entry
\[
\begin{aligned}
	(U_D)_{x,(\xi,j)}
	&=c^{-1/2}\zeta_{p^a}^{x_p y_{p,\xi}}
	\zeta_{q^b}^{x_q(z_\xi+jn)}\\
	&=c^{-1/2}
	\bigl(\zeta_{p^a}^{x_p y_{p,\xi}}\zeta_{q^b}^{x_qz_\xi}\bigr)
	\zeta_{q^b}^{x_qjn}\\
	&=c^{-1/2}\chi_{y_\xi}(x)\zeta_q^{jx_q}
	=c^{-1/2}\chi_{y_\xi}(x)\zeta_q^{j\alpha}.
\end{aligned}
\]
Taking the $s$th Fourier output of these $q$ entries gives
\[
\begin{aligned}
(U_DV)_{x,(s,\xi)}
&=\frac{1}{\sqrt q}\sum_{j=0}^{q-1}
  (U_D)_{x,(\xi,j)}\zeta_q^{-js}\\
&=\frac{\chi_{y_\xi}(x)}{\sqrt{cq}}
  \sum_{j=0}^{q-1}\zeta_q^{j(\alpha-s)}
=\begin{cases}
  \mu^{-1/2}\chi_{y_\xi}(x),&s=\alpha,\\
  0,&s\ne\alpha.
 \end{cases}
\end{aligned}
\]
Thus the row indexed by $x$ is supported in channel $\alpha$.

Let $X=X_{v,\epsilon,\alpha}$, and for $x\in X$ let $r_x$
be the row of $U_DV$ restricted to channel $\alpha$. In
$V^*K_XV=(U_DV)^*E_X(U_DV)$, only the rows indexed by $X$
contribute. Each of their outer products is supported in the
$(\alpha,\alpha)$ block, so the entire matrix has just one
possibly nonzero block
\[
Q_X=\sum_{x\in X}r_x^*r_x .
\]
By Lemma~\ref{prop:row-class-projections}, $K_X$ is an orthogonal
projection. Its unitary conjugate $V^*K_XV$ is therefore an
orthogonal projection, and so is its block $Q_X$.

Since $y_\xi\in D_{v,u}$, its $p$-coordinate has the form
$y_{p,\xi}=u+p^{a-v}n_{y_\xi}$ with $n_{y_\xi}\in\{0,\ldots,p^v-1\}$.
For $x\in X$, use $x_p\equiv\epsilon\pmod{p^v}$ and
$x_q=\alpha+q\,t(x)$ to obtain
\begin{equation}\label{eq:spectralphase}
\begin{aligned}
\chi_{y_\xi}(x)
&=\zeta_{p^a}^{x_p(u+p^{a-v}n_{y_\xi})}
  \zeta_{q^b}^{(\alpha+q\,t(x))z_\xi}\\
&=\zeta_{p^a}^{x_pu}
  \zeta_{p^v}^{\epsilon n_{y_\xi}}
  \zeta_{q^b}^{\alpha z_\xi}
  \zeta_n^{t(x)z_\xi}.
\end{aligned}
\end{equation}
Define the diagonal unitary $\Delta_{\alpha,\epsilon}$ and the
character row $u_t$ on $F_{v,u}$ by
\[
\Delta_{\alpha,\epsilon}
=\diag_{\xi\in I_{v,u}}
   \bigl(\zeta_{p^v}^{\epsilon n_{y_\xi}}\zeta_{q^b}^{\alpha z_\xi}\bigr),
\qquad
u_t=(\zeta_n^{t z_\xi})_{\xi\in I_{v,u}}.
\]
It follows that
\[
r_x=\mu^{-1/2}\zeta_{p^a}^{x_pu}u_{t(x)}\Delta_{\alpha,\epsilon}.
\]
The factor $\zeta_{p^a}^{x_pu}$ cancels in $r_x^*r_x$.
Collecting the terms with the same $t(x)$ gives
\begin{equation}\label{eq:channelgram}
\begin{aligned}
Q_X
&=\Delta_{\alpha,\epsilon}^*
  \left(\frac{1}{\mu}\sum_{x\in X}u_{t(x)}^*u_{t(x)}\right)
  \Delta_{\alpha,\epsilon}\\
&=\Delta_{\alpha,\epsilon}^*
  \left(\sum_{t\in\Z_n}w_X(t)u_t^*u_t\right)\Delta_{\alpha,\epsilon}\\
&=\Delta_{\alpha,\epsilon}^*P(F_{v,u},w_X)\Delta_{\alpha,\epsilon}.
\end{aligned}
\end{equation}
The coefficient in the second line is $w_X(t)$ by
\eqref{eq:spectralweights}. Since $\Delta_{\alpha,\epsilon}$ is
unitary, $P(F_{v,u},w_X)=\Delta_{\alpha,\epsilon}Q_X
\Delta_{\alpha,\epsilon}^*$ is an orthogonal projection.

At $v=0$, we have $X_{0,0,\alpha}=A_\alpha$ and
$\sum_{\alpha=0}^{q-1}K_{A_\alpha}=I_D$ by
Lemma~\ref{prop:row-class-projections}. Using the same $V$ in
each term and comparing channels, we get
\[
\bigoplus_{\alpha=0}^{q-1}Q_{A_\alpha}
=\sum_{\alpha=0}^{q-1}V^*K_{A_\alpha}V
=V^*I_DV
=\bigoplus_{\alpha=0}^{q-1}I_{F_{0,u}} .
\]
Hence $Q_{A_\alpha}=I_{F_{0,u}}$, and \eqref{eq:channelgram}
gives $P(F_{0,u},w_{A_\alpha})=I_{F_{0,u}}$.
\end{proof}

\subsection{Refining one row class}

\begin{lemma}\label{lem:depth}
Let $0\le v\le k$ and suppose $F_{v,u}\ne\varnothing$.
For every $X=X_{v,\epsilon,\alpha}$, the weight $w_X$ in
\eqref{eq:spectralweights} satisfies $\Rp_v$ on $F_{v,u}$.
\end{lemma}

\begin{proof}
We induct on $v$. At $v=0$, the row classes are $A_\alpha$.
Lemma~\ref{lem:grams} gives $P(F_{0,u},w_{A_\alpha})=I_{F_{0,u}}$,
so $w_{A_\alpha}$ satisfies $\Rp_0$ on $F_{0,u}$ by
Lemma~\ref{lem:saturated}.

Suppose $1\le v\le k$ and the assertion holds at level $v-1$.
Fix $\alpha,u$ and put $D=S\cap D_{v,u}$ and $F=F_{v,u}$, with
nonempty index set $I=I_{v,u}$. For $0\le\epsilon<p^{v-1}$, let
\[
X_{\mathrm{par}}=X_{v-1,\epsilon,\alpha},\qquad
\epsilon_i=\epsilon+i p^{v-1},\qquad
X_i=X_{v,\epsilon_i,\alpha}\quad(0\le i<p).
\]
$X_{\mathrm{par}}=\bigsqcup_{i=0}^{p-1}X_i$.
If $X_{\mathrm{par}}=\varnothing$, all the child weights vanish,
and their row graphs satisfy $\Rp_v$ because they are empty.
Assume $X_{\mathrm{par}}\ne\varnothing$. Put $w_i=w_{X_i}$ and
$P_i=P(F,w_i)$. The weights are nonnegative and rational by
\eqref{eq:spectralweights}, and each $P_i$ is an orthogonal
projection by Lemma~\ref{lem:grams}.

For $\xi\in I$, write $y_{p,\xi}=u+p^{a-v}n_{y_\xi}$ as in
Section~4.3. Partition $I$ by
\[
I_\beta=\{\xi\in I:n_{y_\xi}\equiv\beta\pmod p\},
\qquad u_\beta=u+\beta p^{a-v},\qquad 0\le\beta<p.
\]
Since $y_{p,\xi}-u_\beta=p^{a-v}(n_{y_\xi}-\beta)$, the condition
$\xi\in I_\beta$ is equivalent to
$y_{p,\xi}\equiv u_\beta\pmod{p^{a-v+1}}$ for $\xi\in I$.
Take $y\in D$ and let $\xi$ be the index of the cycle containing
$y$. Every point of this cycle has $p$-coordinate $y_{p,\xi}$, so
\[
\xi\in I_\beta
\quad\Longleftrightarrow\quad
y_p=y_{p,\xi}\equiv u_\beta\pmod{p^{a-v+1}}
\quad\Longleftrightarrow\quad
y\in D_{v-1,u_\beta}.
\]
Thus the union of the cycles indexed by $I_\beta$ is
$D\cap D_{v-1,u_\beta}=S\cap D_{v-1,u_\beta}$. 
Here $D_{v-1,u_\beta}\subseteq D_{v,u}$ since
$u_\beta\equiv u\pmod{p^{a-v}}$. Also $0\le u_\beta<p^{a-v+1}$,
so $u_\beta$ is the representative at level $v-1$. With the fixed
cycle order, the sublist indexed by $I_\beta$ is therefore
\[
F_\beta=(z_\xi)_{\xi\in I_\beta}=F_{v-1,u_\beta}.
\]

Let $\omega=\zeta_p$ and set $(\Omega_i)_{\xi,\xi}=\omega^{-i\beta}$
for $\xi\in I_\beta$, as in \eqref{eq:twist}. On the full list $F$,
write $\Delta_i=\Delta_{\alpha,\epsilon_i}$ and let $Q_i$ be the
channel $\alpha$ block of $V^*K_{X_i}V$. The column set $D$ and
the matrix $V$ are common to all the children. By
\eqref{eq:channelgram}, $Q_i=\Delta_i^*P_i\Delta_i$, and
$V^*K_{X_i}V$ has zero entries outside its $(\alpha,\alpha)$ block.
For $i\ne j$, \eqref{eq:gramproduct} gives $K_{X_i}K_{X_j}=0$.
Thus
\[
(V^*K_{X_i}V)(V^*K_{X_j}V)=V^*K_{X_i}K_{X_j}V=0,
\]
whose $(\alpha,\alpha)$ block is $Q_iQ_j=0$.

For $\xi\in I_\beta$, the $q$-coordinate phases cancel in
$\Delta_i\Delta_j^*$. Since $\epsilon_i-\epsilon_j=(i-j)p^{v-1}$,
\[
\begin{aligned}
(\Delta_i\Delta_j^*)_{\xi,\xi}
&=\zeta_{p^v}^{(\epsilon_i-\epsilon_j)n_{y_\xi}}
 =\zeta_p^{(i-j)n_{y_\xi}}\\
&=\omega^{(i-j)\beta}
 =(\Omega_i^*\Omega_j)_{\xi,\xi}.
\end{aligned}
\]
The matrices are diagonal, so $\Delta_i\Delta_j^*=\Omega_i^*\Omega_j$.
Multiplying $Q_iQ_j=0$ on the left by $\Delta_i$ and on the
right by $\Delta_j^*$ gives
\[
0=\Delta_iQ_iQ_j\Delta_j^*
 =P_i\Delta_i\Delta_j^*P_j
 =P_i\Omega_i^*\Omega_jP_j
 \qquad(i\ne j).
\]

\Needspace{11\baselineskip}
The parent and its children use the same $t(x)$ and the same
normalization $\mu$ in \eqref{eq:spectralweights}. For $t\in\Z_n$,
\[
\sum_{i=0}^{p-1}w_i(t)
=\frac1\mu\sum_{i=0}^{p-1}\bigl|\{x\in X_i:t(x)=t\}\bigr|
=w_{X_{\mathrm{par}}}(t).
\]
Put $H_\beta=(\sum_i P_i)[I_\beta,I_\beta]$. For $I_\beta\ne\varnothing$,
the character row on $F_\beta$ is $u_t|_{I_\beta}$. Thus
\[
\begin{aligned}
	H_\beta
	&=\left(\sum_i\sum_{t\in\Z_n}w_i(t)u_t^*u_t\right)
	[I_\beta,I_\beta]\\
	&=\sum_{t\in\Z_n}\left(\sum_i w_i(t)\right)
	(u_t|_{I_\beta})^*(u_t|_{I_\beta})\\
	&=\sum_{t\in\Z_n}w_{X_{\mathrm{par}}}(t)
	(u_t|_{I_\beta})^*(u_t|_{I_\beta})\\
	&=P(F_\beta,w_{X_{\mathrm{par}}})
	=P(F_{v-1,u_\beta},w_{X_{\mathrm{par}}}).
\end{aligned}
\]
Lemma~\ref{lem:grams} gives that this matrix is an orthogonal
projection, and the induction hypothesis gives $\Rp_{v-1}$ for
$w_{X_{\mathrm{par}}}$ on $F_\beta$. An empty part gives the
zero-dimensional projection.

Choose a nonempty part $I_{\beta_0}$, which exists because
$I\ne\varnothing$. Its diagonal entries satisfy
\[
(H_{\beta_0})_{\xi,\xi}
=\sum_{t\in\Z_n}w_{X_{\mathrm{par}}}(t)
=\frac{|X_{\mathrm{par}}|}{\mu}>0
\qquad(\xi\in I_{\beta_0}).
\]
Thus $H_{\beta_0}\ne0$. Applying Theorem~\ref{thm:weights}
with $d=v-1$ shows that every $w_i$ satisfies $\Rp_v$ on
$F_{v,u}$. Varying $\epsilon$ and $i$ covers all row classes
at level $v$, completing the induction.
\end{proof}

\section{Residue proportions and the contradiction}\label{sec:residues}

\subsection{Denominators of residue proportions}

For a nonempty residue class $C$ of $A$ modulo a power of $p$, put
\[
\pi_\alpha(C)=\frac{|C\cap A_\alpha|}{|C|},
\qquad 0\le\alpha<q.
\]

\begin{lemma}\label{lem:proportions}
For $0\le v\le k$ and $C=C_{v,\epsilon}$, every proportion
$\pi_\alpha(C)$ belongs to $\Zp$.
\end{lemma}

\begin{proof}
Fix $C=C_{v,\epsilon}$ and put $X=C\cap A_\alpha$.
By Lemma~\ref{lem:cycles}, $S$ contains a complete $q$-cycle.
Its common $p$-coordinate determines a residue $u$ modulo
$p^{a-v}$, so this cycle lies in $S\cap D_{v,u}$ and
$F_{v,u}\ne\varnothing$.

By Lemma~\ref{lem:depth}, the weight $w_X$ satisfies $\Rp_v$ on
$F_{v,u}$. The components $G$ of its row graph partition
the positive-weight points, so
\[
w_X(\Z_n)=\sum_G w_X(G)\in p^{-v}\Zp,
\qquad n=q^{b-1}.
\]
By \eqref{eq:sizes}, $|C|=c/p^v$ and $\mu=c/q$. Using
$w_X(\Z_n)=|X|/\mu$, we obtain
\begin{equation}\label{eq:pdenom}
\pi_\alpha(C)
=\frac{|X|}{c/p^v}
=\frac{p^v}{q}\frac{|X|}{\mu}
=\frac{p^v}{q}w_X(\Z_n)\in\Zp.
\end{equation}
The membership follows from $p^v w_X(\Z_n)\in\Zp$ and $1/q\in\Zp$.
\end{proof}

Write
\[
\Z[1/p]=\left\{\frac{h}{p^e}:h\in\Z,\ e\in\Z_{\ge0}\right\}.
\]

\begin{lemma}\label{lem:differences}
For $0\le v\le k$, $C=C_{v,\epsilon}$, and $0\le\alpha,\beta<q$,
\[
\pi_\alpha(C)-\pi_\beta(C)\in\Z[1/p].
\]
\end{lemma}

\begin{proof}
Let $K=m+1$ be as in \eqref{eq:firstzeros}, so that
$k<K\le a$, $A[p^K]\ne0$, and $A[p^Kq]=0$.
Extend the row-class definition to depth $K$ by setting
\[
C_{j,r}=\{x\in A:x_p\equiv r\pmod{p^j}\},\qquad
0\le j\le K,\quad 0\le r<p^j.
\]
For $1\le j\le K$, a nonempty parent $C'=C_{j-1,r}$ is the
disjoint union of its children
\[
C'=\bigsqcup_{i=0}^{p-1}C_{j,r+i p^{j-1}} .
\]

\begin{enumerate}[label=\textup{(\alph*)},ref=\textup{(\alph*)}]
\item\label{item:equal}
If $A[p^j]=0$, the $p$ children of $C'$ have equal size.
For $x\in C_{j,s}$, we have $x\equiv x_p\equiv s\pmod{p^j}$.
Grouping the mask terms by these classes gives
\[
0=A[p^j]
 =\sum_{s=0}^{p^j-1}|C_{j,s}|\zeta_{p^j}^s
 =\sum_{r'=0}^{p^{j-1}-1}\sum_{i=0}^{p-1}
   |C_{j,r'+i p^{j-1}}|\,\zeta_{p^j}^{r'+i p^{j-1}} .
\]
By Lemma~\ref{thm:coeff}, the $p$ counts at each $r'$ are equal.
For $r'=r$ these counts sum to $|C'|$, so each child of $C'$
has size $|C'|/p>0$.

\item\label{item:single}
If $A[p^j]\ne0$ and $C'$ meets two or more $q$-residues,
then $C'$ has exactly one nonempty child.
Suppose two children are nonempty and choose two points of
$C'$ with different $q$-residues. If the two points $x_1,x_2$ lie in the same child, choose $z$
from a different nonempty child. Since
\[
(x_1-z)-(x_2-z)=x_1-x_2\not\equiv 0\pmod q,
\]
we have $x_1\not\equiv z\pmod q$ or $x_2\not\equiv z\pmod q$. Thus there are $x,x'\in C'$ in distinct children
with $x\not\equiv x'\pmod q$. For $h=x-x'$ this gives
\[
p^{j-1}\mid h,\qquad p^j\nmid h,\qquad q\nmid h.
\]
Write $h=p^{j-1}h_0$ and $d=N/p^{j-1}$. Then $\gcd(h_0,d)=1$.
Row orthogonality gives
\[
0=\sum_{y\in B}\zeta_N^{hy}
 =\sum_{y\in B}\zeta_d^{h_0y}
 =M_B(\zeta_d^{h_0}).
\]
Since $\zeta_d^{h_0}$ is a primitive $d$th root, $\Phi_d$
divides $M_B$, giving $B[N/p^{j-1}]=B[d]=0$. Since $0\le j-1\le m$,
\eqref{eq:pattern} gives $A[p^j]=0$, a contradiction.
\end{enumerate}

Starting at $C=C_{v,\epsilon}$, stop at a class if it is contained
in a single $A_\gamma$ or has depth $K$. Otherwise, split it into
its nonempty children. When a class $C'$ at depth $j-1$ is split, it meets at least two
residue classes modulo $q$. If $A[p^j]=0$, rule~\ref{item:equal}
gives $p$ nonempty children, each of size $|C'|/p$. If $A[p^j]\ne0$,
rule~\ref{item:single} gives exactly one nonempty child. This child
contains every point of $C'$, so its size is $|C'|$.
Each step increases the depth by one, and the children retained
at each step partition their parent. Thus the set $\mathcal L$
of leaves is finite and
\[
C=\bigsqcup_{D\in\mathcal L}D.
\]
Along the path $C=C^{(0)},C^{(1)},\ldots,C^{(\ell)}=D$ to a leaf,
let $e(D)$ count the steps using rule~\ref{item:equal}. Then
\[
\frac{|D|}{|C|}
=\prod_{i=0}^{\ell-1}\frac{|C^{(i+1)}|}{|C^{(i)}|}
=p^{-e(D)} .
\]
When $C$ is itself a leaf, this is the empty product with
$e(C)=0$.

\Needspace{10\baselineskip}
Put $\rho=\zeta_{p^Kq}$, $\xi=\zeta_{p^K}$, and $E=\Q(\zeta_q)$.
The root $\xi\zeta_q=\rho^{q+p^K}$ is primitive because
$q+p^K$ is coprime to $p^Kq$. Hence $A[p^Kq]=0$ implies
$M_A(\xi\zeta_q)=0$. With
$S_r=\sum_{x\in C_{K,r}}\zeta_q^x\in E$, the congruence
$x\equiv x_p\equiv r\pmod{p^K}$ for $x\in C_{K,r}$ gives
\[
0=M_A(\xi\zeta_q)=\sum_{x\in A}\xi^x\zeta_q^x
 =\sum_{r=0}^{p^K-1}S_r\xi^r.
\]
Since $\gcd(p^K,q)=1$, we have
\[
E(\xi)
=\Q(\zeta_q,\zeta_{p^K})
=\Q(\zeta_{p^Kq})
=\Q(\rho).
\] and
\[
[E(\xi):E]=\frac{\varphi(p^Kq)}{\varphi(q)}=\varphi(p^K).
\]
The minimal polynomial of $\xi$ over $E$ divides $\Phi_{p^K}$
and has the same degree, so it equals $\Phi_{p^K}$.
Hence $\Phi_{p^K}$ is irreducible over $E$.

Writing $r=s+i p^{K-1}$ in $\sum_{r=0}^{p^K-1}S_r\xi^r=0$
and applying Lemma~\ref{thm:coeff} over $E$, we obtain
\[
S_s=S_{s+p^{K-1}}=\cdots=S_{s+(p-1)p^{K-1}},
\qquad 0\le s<p^{K-1}.
\]

\Needspace{7\baselineskip}
Let $D$ be a leaf at depth $K$, with parent $C_{K-1,s}$.
By the stopping rule, the parent meets at least two residue classes
modulo $q$. Since $A[p^K]\ne0$, rule~\ref{item:single} shows that
$D=C_{K,s+i_0 p^{K-1}}$ is its only nonempty child.

For $0\le i<p$ with $i\ne i_0$, the class $C_{K,s+i p^{K-1}}$
is empty, so $S_{s+i p^{K-1}}=0$. The equality of the $p$ coefficients
for this parent gives
\[
0=S_{s+i_0 p^{K-1}}
=\sum_{x\in D}\zeta_q^x
=\sum_{\gamma=0}^{q-1}|D\cap A_\gamma|\zeta_q^\gamma.
\]
By Lemma~\ref{thm:coeff} at $\zeta_q$, the $q$ counts
$|D\cap A_\gamma|$ are equal. Since they sum to $|D|$, we have
\[
|D\cap A_\gamma|=\frac{|D|}{q},
\qquad
\pi_\gamma(D)=\frac1q
\qquad (0\le\gamma<q).
\]

For a leaf stopped before depth $K$, the proportions are $0$
or $1$. On a leaf at depth $K$, both $\pi_\alpha(D)$ and
$\pi_\beta(D)$ equal $1/q$. Thus
$d_D=\pi_\alpha(D)-\pi_\beta(D)\in\{-1,0,1\}$ for every leaf.
Because the leaves partition $C$,
\[
\pi_\alpha(C)
=\sum_{D\in\mathcal L}\frac{|D\cap A_\alpha|}{|C|}
=\sum_{D\in\mathcal L}\frac{|D|}{|C|}\pi_\alpha(D).
\]
Subtracting the corresponding identity for $\beta$ and setting
$e_*=\max_{D\in\mathcal L} e(D)$ gives
\[
\pi_\alpha(C)-\pi_\beta(C)
=\sum_{D\in\mathcal L}p^{-e(D)}d_D
=\frac{\sum_{D\in\mathcal L}d_Dp^{e_*-e(D)}}{p^{e_*}}
\in\Z[1/p].
\qedhere
\]
\end{proof}

\Needspace{10\baselineskip}
\subsection{Proof of Theorem~\texorpdfstring{\ref*{thm:main}}{1}}

\begin{proof}[Proof of Theorem~\ref{thm:main}]
Let $(A,B)$ be the spectral pair chosen in Section~\ref{sec:obstruction}
under the assumption that Theorem~\ref{thm:main} fails.
Fix $C=C_{v,\epsilon}$ with $0\le v\le k$. By
Lemma~\ref{lem:proportions}, both $\pi_\alpha(C)$ and $\pi_\beta(C)$
lie in $\Zp$, so their difference also lies in $\Zp$.
By Lemma~\ref{lem:differences}, the difference also lies in $\Z[1/p]$.
Thus
\[
\pi_\alpha(C)-\pi_\beta(C)\in\Zp\cap\Z[1/p]=\Z
\qquad (0\le\alpha,\beta<q).
\]
The intersection is $\Z$ because a reduced denominator that is both
a power of $p$ and coprime to $p$ must be $1$.

Each proportion $\pi_\gamma(C)$ lies in $[0,1]$.
Since $C=\bigsqcup_{\gamma=0}^{q-1}(C\cap A_\gamma)$, we have
\[
\sum_{\gamma=0}^{q-1}\pi_\gamma(C)
=\frac{1}{|C|}\sum_{\gamma=0}^{q-1}|C\cap A_\gamma|
=\frac{|C|}{|C|}=1.
\]
If all the proportions are equal, then $q\,\pi_\gamma(C)=1$
for every $\gamma$, so $C$ is uniform modulo $q$.

Otherwise, choose $\alpha,\beta$ with
$\pi_\alpha(C)>\pi_\beta(C)$. The difference is a positive
integer at most $1$, so
$\pi_\alpha(C)-\pi_\beta(C)=1$. Since both proportions lie
in $[0,1]$, this forces $\pi_\alpha(C)=1$ and $\pi_\beta(C)=0$.
Thus $|C\cap A_\alpha|=|C|$, so $C\subseteq A_\alpha$.

By \eqref{eq:sizes}, $\pi_\alpha(A)=|A_\alpha|/|A|=1/q$ for every
$\alpha$, so the root class $A=C_{0,0}$ is uniform. Let $1\le v\le k$
and suppose $C=C_{v-1,\epsilon}$ is uniform. Its children are
$C_i=C_{v,\epsilon+i p^{v-1}}$ for $0\le i<p$, and \eqref{eq:sizes}
gives $|C_i|=|C|/p$. Let $t$ be the number of uniform children and
$n_\alpha$ the number contained in $A_\alpha$. For a fixed $\alpha$,
the values $\pi_\alpha(C_i)$ are $1/q$ on the $t$ uniform children,
$1$ on the $n_\alpha$ children contained in $A_\alpha$, and $0$ on
the remaining children. Since the children partition $C$,
\[
\frac1q=\pi_\alpha(C)
=\sum_{i=0}^{p-1}\frac{|C_i|}{|C|}\pi_\alpha(C_i)
=\frac1p\sum_{i=0}^{p-1}\pi_\alpha(C_i)
=\frac1p\left(\frac{t}{q}+n_\alpha\right).
\]
Thus $qn_\alpha=p-t$. Since $0\le p-t\le p<q$ and $q\mid(p-t)$,
we have $t=p$, so all the children are uniform. By induction,
every class at depth $k$ is uniform.

Each class at depth $k$ now meets every residue class modulo $q$.
Since $A[p^{k+1}]\ne0$ by \eqref{eq:firstzeros}, rule~\ref{item:single}
in the proof of Lemma~\ref{lem:differences}, applied with $j=k+1$,
shows that each class at depth $k$ has exactly one nonempty child.
Since $k<a$, we have $p^{k+1}\mid p^a\mid N$. For $v=k,k+1$,
the congruence $x_p\equiv x\pmod{p^v}$ gives
\[
C_{v,r}=\{x\in A:x\equiv r\pmod{p^v}\}
=A\cap(r+p^v\Z_N).
\]
Thus our classes at depths $k$ and $k+1$ are exactly the residue
classes used in property~(iii) of Section~\ref{sec:obstruction}.
Applying  property~(iii) of Section~\ref{sec:obstruction} at $j=k$ gives a residue $\epsilon$ for
which at least two of the classes $C_{k+1,\epsilon+i p^k}$,
$0\le i<p$, are nonempty. This contradicts the preceding
conclusion that $C_{k,\epsilon}$ has exactly one nonempty child
and proves Theorem~\ref{thm:main}.
\end{proof}

\section*{Declaration on the use of generative AI}

OpenAI Codex suggested proof ideas and helped the authors develop proofs
for some of the theorems, locate relevant
publications. The authors take full responsibility for the manuscript, including its
mathematical content and the accuracy of its references.

\medskip
\begingroup
\small
\setlength{\parskip}{0pt}
\noindent\textsc{Zi-Jian Song}\par
\noindent Hubei Key Laboratory of Mathematical Sciences, College of Mathematics
and Statistics, Central China Normal University, Wuhan, Hubei 430079, China\par
\noindent\textit{Email address:} \href{mailto:songzijian2025@outlook.com}{\texttt{songzijian2025@outlook.com}}\par
\smallskip
\noindent\textsc{Wei-Jie Wang}\par
\noindent Hubei Key Laboratory of Mathematical Sciences, College of Mathematics
and Statistics, Central China Normal University, Wuhan, Hubei 430079, China\par
\noindent\textit{Email address:} \href{mailto:wwjmath@163.com}{\texttt{wwjmath@163.com}}\par
\endgroup
\end{document}